\documentclass[12pt]{article}
\usepackage{amsmath}
\usepackage{amssymb}
\numberwithin{equation}{section}
\usepackage{epic}

\title{PDE's for the Gaussian ensemble with external source
and
  the Pearcey distribution}

\author{ Mark Adler\thanks
{2000 {\em Mathematics Subject Classification}. Primary:
60J60, 60J65, 60G55; secondary: 35Q53, 35Q58. {\em Key
words and Phrases}: Non-intersecting Brownian motions,
Pearcey distribution, matrix models, random Hermitian
ensembles, multi-component KP equation, Virasoro
constraints.
\newline
  Department of Mathematics, Brandeis University,
Waltham, MA 02454, USA. E-mail: adler@brandeis.edu. The
support of a National Science Foundation grant \#
DMS-04-06287 is gratefully acknowledged.}~~~~~~ Pierre
van Moerbeke\thanks{ Department of Mathematics,
Universit\'e de Louvain, 1348 Louvain-la-Neuve, Belgium
and Brandeis University, Waltham, MA 02454, USA. E-mail:
vanmoerbeke@math.ucl.ac.be and @brandeis.edu. The
support of a National Science Foundation grant \#
DMS-04-06287, a European Science Foundation grant
(MISGAM), a Marie Curie Grant (ENIGMA), Nato, FNRS and
Francqui Foundation grants is gratefully acknowledged.}}

\date{}

\newcommand{\MAT}[1]{\left(\begin{array}{*#1c}}
\newcommand{\mat}{\end{array}\right)}
\newcommand{\qed}{\leavevmode\unskip\nobreak\penalty200\hskip2pt\null
\nobreak\hfill\rule{1.1ex}{1.1ex}
\medbreak }

\newcommand{\rg}{\rightarrow}

\newcommand{\LR}{{\cal L}}

\newcommand{\BB}{{\cal B}}
\newcommand{\BC}{{\mathbb C}}

\newcommand{\BP}{{\mathbb P}}

\newcommand{\BV}{{\mathbb V}}

\newcommand{\iy}{\infty}
\newcommand{\pl}{\partial}
\newcommand{\al}{\alpha}

\newcommand{\gs}{{\bf s}}
\newcommand{\no}{\nonumber}

\newenvironment{proof}{\medskip\noindent{\it Proof:\/} }{\qed}
\newenvironment
         {remark}{\medskip\noindent\underline{\it Remark:\/} }{\medbreak}

\newcommand{\om}{\omega}
\newcommand{\vp}{\varphi}

\newcommand{\dt}{\delta}
\newcommand{\Dt}{\Delta}
  \newcommand{\vr}{\varepsilon}
\newcommand{\sg}{\sigma}
\newcommand{\BR}{{\mathbb R}}
\newcommand{\lb}{\lambda}
\newcommand{\Lb}{\Lambda}

\newcommand{\BJ}{{\mathbb J}}

\newcommand{\diag}{\operatorname{diag}}

\def\be#1\ee{\begin{equation}#1\end{equation}}
\def\bea#1\eea{\begin{eqnarray}#1\end{eqnarray}}
\def\bean#1\eean{\begin{eqnarray*}#1\end{eqnarray*}}

\newcommand{\Tr}{\operatorname{\rm Tr}}

\newtheorem{definition}{Definition}[section]
\newtheorem{theorem}[definition]{Theorem}
\newtheorem{lemma}[definition]{Lemma}
\newtheorem{corollary}[definition]{Corollary}
\newtheorem{proposition}[definition]{Proposition}

\catcode `!=11

\newdimen\squaresize
\newdimen\thickness
\newdimen\Thickness
\newdimen\ll! \newdimen \uu! \newdimen\dd! \newdimen \rr! \newdimen
\temp!

\def\sq!#1#2#3#4#5{%
\ll!=#1 \uu!=#2 \dd!=#3 \rr!=#4
\setbox0=\hbox{%
  \temp!=\squaresize\advance\temp! by .5\uu!
  \rlap{\kern -.5\ll!
  \vbox{\hrule height \temp! width#1 depth .5\dd!}}%
  \temp!=\squaresize\advance\temp! by -.5\uu!
  \rlap{\raise\temp!
  \vbox{\hrule height #2 width \squaresize}}%
  \rlap{\raise -.5\dd!
  \vbox{\hrule height #3 width \squaresize}}%
  \temp!=\squaresize\advance\temp! by .5\uu!
  \rlap{\kern \squaresize \kern-.5\rr!
  \vbox{\hrule height \temp! width#4 depth .5\dd!}}%
  \rlap{\kern .5\squaresize\raise .5\squaresize
  \vbox to 0pt{\vss\hbox to 0pt{\hss $#5$\hss}\vss}}%
}
  \ht0=0pt \dp0=0pt \box0
}

\def\vsq!#1#2#3#4#5\endvsq!{\vbox to \squaresize{\hrule
width\squaresize height 0pt%
\vss\sq!{#1}{#2}{#3}{#4}{#5}}}

\newdimen \LL! \newdimen \UU! \newdimen \DD! \newdimen \RR!

\def\vvsq!{\futurelet\next\vvvsq!}
\def\vvvsq!{\relax
   \ifx     \next l\LL!=\Thickness \let\continue=\skipnexttoken!
   \else\ifx\next u\UU!=\Thickness \let\continue=\skipnexttoken!
   \else\ifx\next d\DD!=\Thickness \let\continue=\skipnexttoken!
   \else\ifx\next r\RR!=\Thickness \let\continue=\skipnexttoken!
   \else\def\continue{\vsq!\LL!\UU!\DD!\RR!}%
   \fi\fi\fi\fi
   \continue}

\def\skipnexttoken!#1{\vvsq!}

\def\place#1#2#3{\vbox to 0pt{\vss
\rlap{\kern#1\squaresize
   \raise#2\squaresize\hbox{$#3$}}
\vss}}

\catcode `!=12

\begin{document}
\maketitle


\vspace*{-.5cm}

\tableofcontents

\bigbreak

Br\'ezin and Hikami
\cite{Brezin2,Brezin3,Brezin4,Brezin5} have considered a
random Gaussian Hermitian ensemble with external source,
 $$
\frac{1}{Z} e^{-\frac{1}{2}\Tr (M-A)^2}dM,
  $$
where $M$ is random and $A$ is deterministic. Notice
this matrix ensemble, which had come up in the prior
literature \cite{Pastur, Kazakov,Gross}, ceases to be
unitary invariant. The matrix $A$ is chosen so that the
support of the equilibrium measure has a gap, when the
size of the random matrices tends to infinity. Through a
fine tuning of $A$, the gap can be made to close at the
origin. Br\'ezin and Hikami propose a simple model
having this feature, where the matrix $A$ is diagonal,
with two eigenvalues $a$ and $-a$ of equal multiplicity.
Thus, upon letting the size of the matrix go to $\iy$
and after appropriate rescaling, they discover a new
critical distribution, specified by a kernel involving
Pearcey integrals \cite{Pearcey} and having universality
properties.

P. Zinn-Justin \cite{Zinn,Zinn1} establishes the
determinantal form of the correlation functions for the
eigenvalues of the finite model.
Then extending classical connections between random
matrix theory and non-intersecting random paths,
Aptekarev, Bleher and Kuijlaars in \cite{AptBleKui} give
a non-intersecting Brownian motion interpretation of
this Gaussian ensemble with external source. They also
show that multiple orthogonal polynomials are the right
tools for studying this model and its limit (see
\cite{Van Assche,BleKui1,BleKui2,BleKui3})

The present paper studies the Gaussian Hermitian random
matrix ensemble ${\cal H}_n$ with external source $A$,
given by the diagonal matrix (set $n=k_1+k_2$)
 \be
 A:= \left(\begin{array}{cccccc}
   a\\
   &\ddots& & & &{\bf O}\\
   & &a\\
   & & &-a\\
   &{\bf O}& & &\ddots\\
   & & & & &-a
\end{array}\right)\begin{array}{l}
\updownarrow k_1\\
\\
\\
\\
\updownarrow k_2
\end{array},
  \label{0.1}\ee
  and density
  \be
  \frac{1}{Z_n} 
  e^{-\Tr(\frac{1}{2} M^2-AM) }dM.
  \ee
    Given a disjoint union of intervals
 $
 E:= \bigcup^r_{i=1}[b_{2i-1},b_{2i}]\subset \BR
   ,$
%
%
%
define the algebra of differential operators
 \be{\cal
  B}_k=\sum_1^{2r}b_i^{k+1}\frac{\pl}{\pl b_i}.\ee
Consider the following probability:
  \be
\BP_n(a;E):=\BP(~\mbox{all eigenvalues}~\in E) =
 \frac{1}{Z_n} \int_{{\cal H}_{n}(E)}
   e^{-\Tr(\frac{1}{2} M^2-AM) }dM
  \label{0.4}
   \ee
 In \cite{AvM0}, we have shown that for $A=0$,
 the probability for this Gaussian Hermitian ensemble (GUE) satisfies
 a {\em fourth-order PDE, with quadratic non-linearity}
 (reducing to Painlev\'e IV in the case of one boundary point)
 :
  $$ \Bigl({\cal B} _{-1}^4+(4n+6{\cal
B}_{-1}^2\log \BP_n){\cal B}_{-1}^2+ 3{\cal B}^2_0-4
 {\cal B}_{-1}{\cal B}_1+6{\cal B}_0 \Bigr)\log \BP_n=0
 .$$
The {\bf first question} in this paper: {\em Does the
integral (\ref{0.4}), with $A$ as in (\ref{0.1}),
satisfy a PDE?}~ Indeed, we prove:

%
%
%

\begin{theorem} 
%
%
%
  The log of the probability $\BP_n(a;E)$
   satisfies
  a  \underline {fourth-order} PDE in $a$ and in the endpoints $b_1,..., b_{2r}$
  of the set $E$, with \underline {quartic non-linearity}:
  \bea  &&
\Bigl(F^+\BB_{-1}G^-+F^-\BB_{-1}G^+ \Bigr)
\Bigl(F^+\BB_{-1}F^- -F^-\BB_{-1}F^+ \Bigr) \no\\
  &&-
  \Bigl(F^+ G^- +F^- G^+
\Bigr) \Bigl(F^+\BB_{-1}^2F^- -F^-\BB_{-1}^2F^+ \Bigr)
=0, \label{0.5}\eea
%
where\footnote{in terms of the Wronskians
$\{f,g\}_X=gXf-fXg$.}
  \bean
      F^+&:=&  2\BB_{-1}(\frac{\pl}{\pl
   a} - \BB_{-1})\log \BP_n-4k_1 ,~~~~~~~~~~
   ~F^-=  F ^+\Bigr|_{\begin{array}{l}
              a\rightarrow -a \\
               k_1 \leftrightarrow k_2
              \end{array}}
   \no\\
    2 G^+ &:=&
    \left\{
      H_1^+,
          F^+
       \right\}_{\BB_{-1}}
       -\left\{
      H^+_2,  F^+
     \right\}_{\pl/\pl a},~~~~~~~
     ~ G^-=  G^+\Bigr|_{\begin{array}{l}
              a\rightarrow -a \\
               k_1 \leftrightarrow k_2
              \end{array}},
    \eean
    with
    \bea
         H_1^+&:=&\!\!\!\! \frac{\pl}{\pl a}
      \left( \BB_{0} -a\frac{\pl}{\pl a}-a \BB_{-1}
      \right)\log \BP_{n}
       +\left(\BB_{0}\BB_{-1}+4\frac{\pl}{\pl a}\right)
       \log \BP_{n}
     \no\\
     &&~~~~~~~~~~~~~~~~~~~~~~~~~~~~~~~~
     ~~~~~~~~~~~~~~~~~~~~~+4k_1\left(a+\frac{k_2}{a}\right)
     \no\\   \no\\ \hspace{-1cm}
       H^+_2
      &:=&\!\!\!\!\frac{\pl}{\pl a}
      \left( \BB_{0} -a\frac{\pl}{\pl a}-a \BB_{-1}
      \right)\log \BP_n
        -\left( \BB_{0} -2a\BB_{-1}
      -2 \right)\BB_{-1}\log \BP_n.
  \no\\ \label{0.7}\eea
\end{theorem}

%
%
%
%
%
%

 \remark: The change of variables $a\mapsto
-a,~k_1 \leftrightarrow k_2$ in the definition of $F^-$
and $G^-$ act at the level of the operators. In fact,
later, it will be clear that $\BP_n( a;E)$ is invariant
under that change of variables.

\medbreak

Again here we provide a natural integrable deformation
of (\ref{0.4}) (section 1). As is well known, the
probability for $A=0$ relates to the standard Toda
lattice and the {\em one-component KP} equation (see
\cite{AvM0}), the spectrum of coupled random matrices to
the 2-Toda lattice and the {\em two-component KP} (see
\cite{AvM1}), whereas we show
that the model (\ref{0.4}) 
relates to the {\em three-component KP}
equation (section 2). This deformation enjoys Virasoro
constraints as well (section 3), which together with the
bilinear relations arising from 3-KP leads to the PDE of
Theorem 0.1 (section 4).


\bigbreak

The {\bf second question} concerns the {\em Pearcey
distribution}, which we now explain.
Following \cite{AptBleKui}, consider $n=2k$
non-intersecting Brownian motions on $\BR$ (Dyson's
Brownian motions),
   all starting at
the origin, such that the $k$ left paths end up at $-a$
and the $k$ right paths end up at $+a$ at time $t=1$.
 As observed in \cite{AptBleKui}, the
Karlin-McGregor formula \cite{Karlin} enables one to
express the transition probability in terms of the
Gaussian Hermitian random matrix probability $\BP(a;E)$
with external source, as in (\ref{0.4}),
\bea   
\BP_0^{\pm a} \left( \mbox{all
$x_j(t)\in E$} \right)
 &=&\lim_{\tiny\begin{array}{c}
  \mbox{\tiny   all $~\gamma_i\rightarrow 0$}
  \\
 \mbox{\tiny $\delta_1, \ldots , \delta_k\rightarrow -a$}\\
  \mbox{\tiny $\delta_{k+1}, \ldots , \delta_{2k}
  \rightarrow a$}
  \end{array}} \int_{E^n}
   \no\\ &&~~~\hspace*{-3cm}
 \frac{1}{Z_{n}}\det(p(t;\gamma_i,x_j))_{1\leq i,j\leq
n}\det(p(1-t;x_{i'},\delta_{j'}))_{1\leq i',j'\leq n}
  \prod_1^n dx_i,  \no\\
 &=& \BP_{n}\left(a \sqrt{\frac{2t }{1-t}};
E\sqrt{\frac{2}{t(1-t)}}\right)\eea
 where $p(t,x,y)$ is the Brownian transition probability
\be
  p(t,x,y):= \frac{1}{\sqrt{ \pi t}
  }~e^{-\frac{(y-x)^2}{ t}}.
  \label{Brownian}\ee
Let now the number $n=2k$ of particles go to infinity,
and let the points $a$ and $-a$ go to $\pm \iy$. This
forces the left $k$ particles to $-\iy$ at $t=1$ and the
right $k$ particles to $+\iy$ at $t=1$. Since the
particles all leave from the origin at $t=0$, it is
natural to believe that for small times the equilibrium
measure (mean density of particles) is supported by one
interval, and for times close to $1$, the equilibrium
measure is supported by two intervals. With a precise
scaling, $t=1/2$ is critical in the sense that for
$t<1/2$, the equilibrium measure for the particles is
supported by one, and for $t>1/2$, it is supported by
two intervals.
 The Pearcey process ${\cal P}(s)$ is now
defined as the motion of an infinite number of
non-intersecting Brownian paths, just after time
$t=1/2$, with the precise scaling (see
\cite{AptBleKui}):
 \be
 n=2k=\frac{2}{z^4},\quad \pm a=\pm\frac{ 1}{z^2},\quad
 x_i\mapsto x_iz,\quad
 t\mapsto\frac{1}{2}+tz^2,~~~\mbox{~for~}z\rg 0.
 \label{0.10}
 \ee
 Even though the pathwise interpretation
of ${\cal P}(t)$ is unclear and still deserves
investigation, it is natural to define the probability
$$
\BP({\cal P}(t) \cap E= \emptyset)
 := \lim_{z\rightarrow
0} \BP_0^{\pm 1/z^2}\left.\left(\mbox{all}~
x_j\bigl(\frac{1}{2}+tz^2\bigr)
\notin  zE
 ; ~1\leq j\leq
n
  \right)\right|_{n=\frac{2}{z^4}}
  . $$
Br\'ezin and Hikami
\cite{Brezin2,Brezin3,Brezin4,Brezin5} for the Pearcey
kernel and Tracy-Widom \cite{TW} for the extended
kernels show that this limit exists and equals a
Fredholm determinant:
$$
 \BP({\cal P}(t) \cap E= \emptyset)
 =
 \det \left( I-K_{t}\chi_{_{E}}
\right),
$$
where $K_t(x,y)$ is the Pearcey kernel, defined as
follows: \bea
K_t(x,y)
&:=&\frac{p(x)q''(y)
 -p'(x)q'(y)+p''(x)q(y)-t p(x)q(y)}{x-y}
 \no\\
 &=& \int^{\iy}_0p(x+z)q(y+z)dz,
\label{0.11}\eea
where (note $\om =e^{i\pi /4}$)
\bean p(x)&:=&\frac{1}{2\pi}\int^{\iy}_{-\iy}
e^{-\frac{u^4}{4}-\frac{tu^2}{2}-iux}du
 \\
q(y)&:=&\frac{1}{2\pi i}\int_X
 e^{ \frac{u^4}{4}-\frac{tu^2}{2}+uy}du
 ={\rm Im}\left[ \frac{\om}{\pi}\int_0^{\iy}du
e^{-\frac{u^4}{4}
 -\frac{it}{2}u^2}(e^{\om
uy}-e^{-\om uy })\right]
  \eean
satisfy the differential equations
%
$$
p'''-tp'-xp=0\mbox{~and~}q'''-tq'+yq=0.
$$
The contour $X$ is given by the ingoing rays from $\pm
\iy e^{i\pi/4}$ to $0$ and the outgoing rays from $0$ to
 $\pm
\iy e^{-i\pi/4}$, i.e., $X$ stands for the contour
 $$ \nwarrow ~  \swarrow 
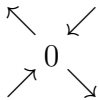
$$

\vspace{-.9cm}

$$ 0
$$

\vspace{-.9cm}

$$  \nearrow~ \searrow   $$
The second result of this paper\footnote{Tracy and Widom
show in \cite{TW} the existence of a large system of
PDE's involving a large system of auxiliary variables
for Q and also for the joint probabilities at different
times.} is to give a PDE for the Pearcey distribution
below in terms of the parameter $t$ appearing in the
kernel (\ref{0.11}). Since this Pearcey distribution
with the parameter $t$ can also be interpreted as the
transition probability for the Pearcey process, we
prove:

\begin{theorem}

For compact $E =\bigcup^r_{i=1}[x_{2i-1},x_{2i}]$ and $
  \BB_j=\sum_1^{2r} x_i^{j+1} \frac{\pl}{\pl x_i}
  $,
%
  \bea
   Q(t;x_1,\ldots,x_{2r})&=&\log \BP
\Bigl(
     {\cal P}(t)\cap E=\emptyset \Bigr)
     =
     \log\det \left(I-K_t\chi_E\right)
  \eea
   satisfies a 4th order
and 3rd degree PDE, which can be written as a
Wronskian\footnote{given that $\{f,g\}_X:=Xf.g-f.Xg$.}:
%
%
   $$
 \left\{  \frac{1}{2}
  \frac{\pl^3 Q}{\pl t^3}+
(\BB_0-2)\BB_{-1}^2Q
  +\frac{1}{16}\Bigl\{
\BB_{-1}\frac{\pl Q}{\pl
t},\BB_{-1}^2Q\Bigr\}_{\BB_{-1}}
  ~,~
   \BB_{-1}^2\frac{\pl Q}{\pl
t}  \right\}_{\BB_{-1}}=0.
$$


\end{theorem}

The proof of this statement, based on taking a limit on
the PDE of Theorem 0.1 will be given in section 5.


\section{An integrable deformation of Gaussian random ensemble with external source}

Consider an ensemble of $n\times n$ Hermitian matrices
with an external source, given by a diagonal matrix
 $$
 A=\diag(a_1,\ldots,a_n)
  $$
 and a general potential $V(z)$, with density
  $$
  \BP_n(M\in [M,M+dM])=
  \frac{1}{Z_{n}}
  e^{-\Tr(V(M)-AM) }dM.
  $$
For the disjoint union of intervals
  $
  E =\bigcup^r_{i=1}[b_{2i-1},b_{2i}]
  $,
the following probability can be transformed by the
Harrish-Chandra-Itzykson-Zuber formula, with
$D=\diag(\lb_1,\ldots,\lb_n$),
 \bea
  \BP_n(\mbox{spectrum} ~M\subset E)&=&\frac{1}{Z_k} \int_{{\cal
  H}_{n}(E)} e^{-\Tr(V(M)-AM) }dM
  \no\\
   && \no\\
  &=&\frac{1}{Z_n} \int_{E^n} \Dt^2_n(\lb)\prod_1^n
e^{-V(\lb_i)}d\lb_i\int_{U(n) }e^{\Tr AUD
U^{-1}}dU  \no\\
&=&\frac{1}{Z'_n} \int_{E^n} \Dt^2_n(\lb)\prod_1^n
e^{-V(\lb_i)}d\lb_i\frac{\det[e^{a_i\lb_j}]_{1\leq i,j\leq
n}}{\Dt_n(\lb)\Dt_n(a)} \no\\
&=&\frac{1}{Z''_n}
\int_{E^n}\Dt_n(\lb)\det[e^{-V(\lb_j)+a_i\lb_j}]_{1\leq
i,j\leq n}\prod_1^n d\lb_i.\no\\
 \label{1.1}\eea
%
  In the following Proposition, we
consider a general situation, of which (\ref{1.1}) with
$A=\diag (a,\ldots,a,-a,\ldots,-a)$ is a special case,
by setting $\varphi^+=e^{az}$ and $\varphi^-=e^{-az}$.
 Consider the Vandermonde determinant of the variables $
 x_1,\ldots,x_{k_1}
 ,y_1,\ldots, y_{k_2}$, namely
 \be\Dt_n(x,y)
 :=\Dt_n(x_1,\ldots,x_{k_1}
 ,y_1,\ldots, y_{k_2}). \label{Vandermonde}\ee
Then we have

\begin{proposition}
Given an arbitrary potential $V(z)$ and arbitrary
functions
 $\varphi^+(z)$ and $\varphi^-(z)$, define
 ($n=k_1+k_2$)
\bean
(\rho_1,\ldots,\rho_n)&:=&e^{-V(z)}\left(\varphi^+(z),z\varphi^+(z)
,\ldots,z^{k_1-1}\varphi^+(z),\right.\\
&&~~~~~~~~~~~~~~~~~~~~~~~ \left.~\varphi^-(z),
z\varphi^-(z),\ldots,z^{k_2-1}\varphi^-(z)\right)
 .
 \eean
 we have

 \bea
 \lefteqn{\frac{1}{n!}\int_{E^n}\Dt_n(z) \det(\rho_i(z_j))_{1\leq i,j\leq
n}
   \prod_1^n
    dz_i}
    \no\\
&=&\frac{ 1}{k_1!k_2!}\int_{E^n}\Dt_n(x,y)
\Dt_{k_1}(x)\Dt_{k_2}(y)\prod_1^{k_1}\varphi^+(x_i)
e^{-V(x_i)}dx_i
\prod_1^{k_2}\varphi^-(y_i)e^{-V(y_i)}dy_i
 \no\\
&=& \det\left(\begin{array}{c}
\left(\displaystyle{\int_{E}}z^{i+j-1}\vp^+(z)e^{-V(z)}\right)_{\tiny{\begin{array}{l}
     1\leq i\leq k_1\\
     0\leq j\leq k_1+k_2-1
     \end{array}}} \\
  \\
\left(\displaystyle{\int_{E}}z^{i+j-1}\vp^-(z)e^{-V(z)}\right)_{\tiny{\begin{array}{l}
     1\leq i\leq k_2\\
     0\leq j\leq k_1+k_2-1
     \end{array}}}
\end{array}\right)
 \label{1.3}
 \eea

\end{proposition}


\proof On the one hand, using
$$
 \det \left( a_{ik}\right)_{1\leq i,k\leq n}
  \det \left( b_{ik}\right)_{1\leq i,k\leq n}=
   \sum_{\sigma \in S_n}\det \left(
a_{i,\sigma(j)}~b_{j,\sigma(j)}\right)_{1\leq i,j\leq n}
,
  $$
and distributing the integration over the different
columns, one computes
  \bea
\lefteqn{\int\Dt_n(z)\det(\rho_i(z_j))_{1\leq i,j\leq n}
 \prod_1^ndz_i}\no\\
&=&\int_{E^n}\det(z_j^{i-1})_{1\leq i,j\leq n}
\det\left(
\begin{array}{c}
     \left(z_j^{i-1}\vp^+(z_j)e^{-V(z_j)}\right)_{\tiny{\begin{array}{l}
     1\leq i\leq k_1\\
     1\leq j\leq k_1+k_2
     \end{array}}}\no\\
     \no\\
   \left(z_j^{i-1}\vp^-(z_j)e^{-V(z_j)}\right)_{\tiny{\begin{array}{l}
     1\leq i\leq k_2\\
     1\leq j\leq k_1+k_2
     \end{array}}}
\end{array}
\right)\prod_1^ndz_i  \no\\
 \no\\
&=&n!\det\left(\begin{array}{c}
\left(\displaystyle{\int_{E}}z^{i+j-1}\vp^+(z)e^{-V(z)}dz\right)_{\tiny{\begin{array}{l}
     1\leq i\leq k_1\\
     0\leq j\leq k_1+k_2-1
     \end{array}}}\no\\
  \\
\left(\displaystyle{\int_{E}}z^{i+j-1}\vp^-(z)e^{-V(z)}dz\right)_{\tiny{\begin{array}{l}
     1\leq i\leq k_2\\
     0\leq j\leq k_1+k_2-1
     \end{array}}}
\end{array}\right)\no\\
  \label{1.4}\eea
On the other hand, one computes
  \bean
\lefteqn{\int_{E^n}\Dt_n(z) \det(\rho_i(z_j))_{1\leq
i,j\leq n}
   \prod_1^n
    dz_i}~\hspace*{15cm}~
  \eean\bean\\
\\
&=& \int_{E^n}\Dt_n(z_1,\ldots,z_n)\sum_{\sg\in
S_n}(-1)^{\sg}\prod^n_{i=1}\rho_i(z_{\sg(i)})
\prod_1^n dz_{\sigma(i)}\\
\\
&=& \sum_{\sigma}(-1)^{\sg}\int_{E^n}
\Dt_n(z_{\sg^{-1}(1)},\ldots,z_{\sg^{-1}(n)})
  \prod_{i=1}^n\rho_i(z_i)
\prod^{n}_1  dz_i\\
\\
&=& \sum_{\sg}(-1)^{\sg}\int_{E^n}(-1)^{\sg}
   \Dt_n(z_1,\ldots,z_n)
  \prod_{i=1}^n\rho_i(z_i)dz_i
   \\
    \\
&=& n!\int_{E^n}\Dt_n(x,y)\prod_1^{k_1}x_i^{i-1}
\varphi^+(x_i)e^{-V(x_i) }dx_i\prod_1^{k_2}y_i^{i-1}
\varphi^-(y_i)e^{-V(y_i)  }dy_i\\
\\
&=&\frac{ n!}{k_1!k_2!}\int_{E^n}\Dt_n(x,y)
\Dt_{k_1}(x)\Dt_{k_2}(y)\prod_1^{k_1}\varphi^+(x_i)
e^{-V(x_i)}dx_i
\prod_1^{k_2}\varphi^-(y_i)e^{-V(y_i)}dy_i,\eean
 where $\Dt_n(x,y)$ is defined in (\ref{Vandermonde}).
 In the last identity, one uses twice the following general identity for a
skew-symmetric function $F(x_1,\ldots,x_{k})$ and a
general measure $\mu(dx)$,
\bean
 \lefteqn{ \int_{\BR^{k }}F(x_1,\ldots, x_{k })\Dt_{k
}(x)\prod_1^{k }\mu(dx_i)
 }\\
 \\
&=&\int_{\BR^{k }}F(x_1,\ldots, x_{k})\sum_{\sg\in
S_{k}}(-1)^{\sg}\prod_{i=1}^{k }x_{\sg(i)}^{i-1}
\prod_1^{k }\mu(
dx_i)\\
\\
&=&\int_{\BR^{k}}\sum_{\sg\in
S_{k}}(-1)^{\sg}F(x_1,\ldots,x_k)
\prod_{i=1}^{k}x_{\sg(i)}^{i-1}
\mu(dx_{\sg(i)})\\
\\
&=&\int_{\BR^{k }}\sum_{\sg\in S_{k }}
(-1)^{\sg}F(x_{\sg^{-1}(1)},\ldots,x_{\sg^{-1}(k)})\prod_{i=1}^{k_1}x_i^{i-1}
\mu(dx_i)\\
&=&k!\int_{\BR^{k }}~F(x_1,\ldots,x_{k
})\prod_{i=1}^{k}x_i^{i-1}\mu(dx_i). \eean
 This ends the proof of Proposition 1.1.\qed


Add extra variables in the exponentials,
  one set for each Vandermonde determinant:
   $$t=(t_1,t_2,\ldots),~~
                       s=(s_1,s_2,\ldots),~~
                       u=(u_1,u_2,\ldots)~~\mbox{and}~~
                       \beta$$
%
%
%
Then, setting (
$n=k_1+k_2$),
 $$
  V(z):=\frac{z^2}{2}+\sum_1^{\iy}t_i z^i ,~~~\varphi^+(z)= e^{ az+\beta z^2
  -\sum_1^{\iy}s_iz^i},
   ~~~\varphi^-(z)= e^{-az-\beta z^2
  -\sum_1^{\iy}u_iz^i},
  $$
   Proposition 1.1 implies
   \bea
    \tau_{k_1k_2}(t,s,u;\beta;E)&:=&
   \det m_{k_1,k_2}(t,s,u;\beta;E) \no\\
   \no\\
   &=&\frac{1}{k_1!k_2!}
  \int_{E^{n }}\Dt_{n }(x,y)
  \prod_{j=1}^{k_1}e^{\sum_1^{\iy}t_i x_j^i  }
   \prod_{j=1}^{k_2}e^{\sum_1^{\iy}t_i y^i_j  }
  \no\\
& &~~~  \left(\Dt_{k_1}(x)\prod^{k_1}_{j=1}
  e^{-\frac{x_j^2}{2}+ax_j+\beta
x_j^2 }
e^{-\sum_1^{\iy}s_ix^i_j}dx_j\right)
 \no\\ &&~~~
\left(\Dt_{k_2}(y)    \prod^{k_2}_{j=1}
e^{-\frac{y_j^2}{2}-ay_j-\beta y_j^2}
e^{-\sum_1^{\iy}u_iy^i_j}dy_j\right)
  ,~~~~~\label{1.5}\eea
where
 \bean
  m_{k_1,k_2}(t,s,u;\beta;E)&:=&
 %
\left(\begin{array}{c}
 \left(\mu^+_{ij}(t,s;\beta,E)\right)_{\tiny{\begin{array}{l}
     1\leq i\leq k_1\\
     0\leq j\leq k_1+k_2-1
     \end{array}}}
     \\
     \left(\mu^-_{ij}(t,u;\beta,E)\right)_{\tiny{\begin{array}{l}
      1\leq i\leq k_2\\
      0\leq j\leq k_1+k_2-1
      \end{array}}}
      \end{array}
      \right),
      \eean
with
  \bea
  \mu^+_{ij}(t,s;\beta,E)&=&
   \displaystyle{\int_{E}}z^{i+j-1}e^{-\frac{z^2}{2}+az+\beta
z^2}e^{\sum_1^{\iy}(t_k-s_k)z^k}dz
  \no \\
  \mu^-_{ij}(t,u;\beta,E)&=&
   \displaystyle{\int_{E}}z^{i+j-1}e^{-\frac{z^2}{2}-az-\beta
z^2} e^{\sum_1^{\iy}(t_k-u_k)z^k}dz
  .\label{1.6}\eea

In particular, by (\ref{1.3}), the integral in
(\ref{1.1}) has the following determinantal
representation in terms of moments:
\bea \lefteqn{\frac{1}{n!}
\int_{E^n}\Dt_n(\lb)\det[e^{-V(\lb_j)+a_i\lb_j}]_{1\leq
i,j\leq n}\prod_1^n d\lb_i}
 \no\\
 &=&\det\left(\begin{array}{c}
\left(\displaystyle{\int_{E}}z^{i+j-1}e^{-\frac{z^2}{2}+az
  } dz\right)_{\tiny{\begin{array}{l}
     1\leq i\leq k_1\\
     0\leq j\leq k_1+k_2-1
     \end{array}}}\\
\\
\left(\displaystyle{\int_{E}}z^{i+j-1}e^{-\frac{z^2}{2}-az
  } dz\right)_{\tiny{\begin{array}{l}
     1\leq i\leq k_2\\
     0\leq j\leq k_1+k_2-1
     \end{array}}}
\end{array}\right)
\label{1.7}\eea  and so
 \be
   \BP_n(\mbox{spec} ~M\subset E)
  =\left.\frac{\tau_{k_1k_2}(t,s,u;\beta;E)}
       {\tau_{k_1k_2}(t,s,u;\beta;\BR)}
       \right|_{t=s=u=\beta=0}
   .\label{1.8}\ee

\remark The integral enjoys the obvious duality:
\be
 x\longleftrightarrow y,~k_1 \longleftrightarrow
k_2,~t\longleftrightarrow t, ~s\longleftrightarrow
u,~a\longleftrightarrow -a,~\beta \longleftrightarrow
-\beta. \label{1.9}\ee


\section{Integrable deformations and 3-component KP
}

\begin{theorem}{\em (Adler-van Moerbeke \cite{AvM4})}
Given the functions $\tau_{n_1,n_2}$ as in (\ref{1.5}),
the wave matrix
 $$
W_{n_1,n_2}^{\pm}(\lb;t,s,u):=\left(\begin{array}{lll}
\Psi_{n_1,n_2}^{(1)\pm}&\Psi_{n_1\pm 1,n_2}^{(2)\pm}&
        \Psi_{n_1,n_2\pm 1}^{(3)\pm}\\
\\
\Psi_{n_1\mp 1,n_2}^{(1)\pm}&\Psi_{n_1,n_2}^{(2)\pm}&
        \Psi_{n_1\mp 1,n_2\pm 1}^{(3)\pm}\\
\\
\Psi_{n_1,n_2\mp 1}^{(1)\pm}&\Psi_{n_1\pm 1,n_2\mp
1}^{(2)\pm}&\Psi_{n_1,n_2}^{(3)\pm}
\end{array}\right)
 $$
 with functions
   \bea
\Psi^{(1)\pm}_{n_1,n_2}(\lb;t,s,u)&:=& \lb^{\pm
(n_1+n_2)}e^{\pm\sum t_i \lb^i}
  \frac{\tau_{n_1,n_2}(t\mp\left[\lb^{-1}\right],s,u)}
      {\tau_{n_1,n_2}(t,s,u)}
     \no \\
  \Psi^{(2)\pm}_{n_1,n_2}(\lb;t,s,u)&:=&
\lb^{\mp n_1}e^{\pm \sum s_i \lb^i}
  \frac{\tau_{n_1,n_2}(t,s\mp \left[\lb^{-1}\right],u)}
      {\tau_{n_1,n_2}(t,s,u)}
     \no  \\
  \Psi^{(3)\pm}_{n_1,n_2}(\lb;t,s,u)&:=&
\lb^{\mp n_2}e^{\pm \sum u_i \lb^i}
  \frac{\tau_{n_1,n_2}(t,s,u\mp \left[\lb^{-1}\right])}
      {\tau_{n_1,n_2}(t,s,u)}
     \eea
     satisfies the bilinear identity
          \be
            \oint_{\iy}
W^+_{k_1,k_2 }(\lb;t,s,u)W_{\ell_1,\ell_2 }^{- }
(\lb;t',s',u')^{\top}d\lb=0
 \label{Psi-bilinear identity}
  \ee
for all integers $k_1,k_2,\ell_1,\ell_2 \geq 0$ and
$t,s,u,t',s',u' \in \BC^{\iy}$.

\end{theorem}

\proof The moments, as defined in (\ref{1.6}), satisfy
$$
\frac{\pl \mu^{\pm}_{ij}}{\pl t_k}=\mu^{\pm}_{i,j+k}
$$
      \bea
\frac{\pl \mu^+_{ij}}{\pl s_k}&=&-\mu^+_{i+k,j},~~~~~~~~
\frac{\pl \mu^+_{ij}}{\pl u_k}=0 \no\\
\frac{\pl \mu^-_{ij}}{\pl s_k}&=&0,~~~~~~~~~~~~~~~~
\frac{\pl \mu^-_{ij}}{\pl u_k}=-\mu^-_{i+k,j},
  \label{2.3}\eea
 and in matrix notation
\be
   \frac{\pl m_{\iy,\iy}}{\pl t_k}=
    m_{\iy,\iy}\Lb^{\top k}
    ,~~
    \frac{\pl m_{\iy,\iy}}{\pl s_k}=
    -\Lb_{-}^{  k} m_{\iy,\iy}
    ,~~
     \frac{\pl m_{\iy,\iy}}{\pl u_k}=
    -\Lb_{+}^{  k} m_{\iy,\iy},
   \label{2.4} \ee
 where
$$\Lb := \left(\begin{tabular}{lllll}
  0& 1&0 &0 \\
  0&0
&1&0 \\
0 &0 &0&1 \\
0 &0 &0&0&  \\
  & & & &$\ddots$
\end{tabular}
\right),~~
  \Lb_-=\left(\begin{tabular}{ll}
   $\Lb$& O\\
   O  & O
   \end{tabular}
   \right),~~
   \Lb_+=\left(\begin{tabular}{ll}
   O& O\\
   O  & $\Lb$
   \end{tabular}
   \right).
$$
Thus the moment matrix satisfies
\bean m_{\iy,\iy}(t,s,u)
&=&e^{-\sum_1^{\iy}(s_k\Lb_-^{
k}+u_k\Lb_+^{k})}m_{\iy,\iy}(0,0,0)e^{\sum_1^{\iy}
t_k\Lb^{\top k}}, \eean
 which implies the
 bilinear identity (\ref{Psi-bilinear identity}).
 The details of proof can be found in \cite{AvM4}.\qed


\begin{corollary} Given the above $\tau$-functions
$\tau_{k_1k_2}(t,s,u)$, they satisfy the bilinear
identities \bean 0\!\!&=&\!\!\!\!\oint_{\iy} d\lb
~\lb^{k_1+k_2-\ell_1-\ell_2}
   e^{ \sum (t_i-t'_i) \lb^i}
    \tau_{k_1,k_2}\left(t-\!\!\left[\lb^{-1}\right],s,u\right)
    \tau_{\ell_1,\ell_2}\left(t'+\!\!\left[\lb^{-1}\right],s',u'\right)
\\\!\!&+ &\!\!\!\!
   \oint d\lb ~\lb^{\ell_1-k_1-2}
   e^{ \sum (s_i-s'_i) \lb^i}
    \tau_{k_1+1,k_2 }\left(t,s-\!\!\left[\lb^{-1}\right] ,u\right)
    \tau_{\ell_1-1 ,\ell_2 }\left(t',s'+\!\!\left[\lb^{-1}\right] ,u'\right)
\\\!\!& +&\!\!\!\!
   \oint d\lb ~\lb^{\ell_2-k_2-2}
   e^{ \sum (u_i-u'_i) \lb^i}
    \tau_{k_1 ,k_2 +1}\left(t,s,u-\!\!\left[\lb^{-1}\right]  \right)
    \tau_{\ell_1 ,\ell_2 -1 }
    \left(t',s',u'+\!\! \left[\lb^{-1}\right]
   \right).
\eean \be\label{tau-bilinear identity} \ee
  Upon specializing, these identities imply PDE's expressed in
  terms of Hirota's symbol~\footnote{Given a polynomial
  $p(t_1,t_2,...)$, define the
customary Hirota symbol $p(\pl_t)f\circ g:=
p(\frac{\pl}{\pl y_1},\frac{\pl}{\pl
y_2},...)f(t+y)g(t-y) \Bigl|_{y=0}$. For later use, the
${\bf s}_{\ell}$'s are the elementary Schur polynomials
$e^{\sum^{\iy}_{1}t_iz^i}:=\sum_{i\geq 0} {\bf
s}_i(t)z^i$ and set ${\bf s}_{\ell}(\tilde \pl):={\bf
s}_{\ell}(\frac{\pl}{\pl t_1},\frac{1}{2}\frac{\pl}{\pl
t_2},\ldots).$ }, for $j=1,2,\ldots$:
 \bea
\gs_j(\tilde\pl_t)\tau_{k_1+1,k_2}\circ\tau_{k_1-1,k_2}&=&\tau^2_{k_1k_2}
\frac{\pl^2}{\pl  s_1\pl t_{j+1}}\log\tau_{k_1k_2}
  \label{bilinear identity1}\\
\no\\
\gs_j(\tilde\pl_s)\tau_{k_1-1,k_2}\circ\tau_{k_1+1,k_2}&=&\tau^2_{k_1k_2}
\frac{\pl^2}{\pl t_1\pl s_{j+1}}\log\tau_{k_1k_2}
  \label{bilinear identity2}
  , \eea
yielding
 \bea \frac{\pl^2\log \tau_{k_1 ,k_2 }}{\pl
t_1\pl s_1} &=&
   \frac{\tau_{k_1 +1,k_2 }\tau_{k_1 -1,k_2 }}
         {\tau_{k_1 ,k_2 }^2}
\label{14}\\
\frac{\pl}{\pl t_1}\log \frac{\tau_{k_1 +1,k_2 }}
                           {\tau_{k_1 -1,k_2 }}
                           &=&\frac
      {\frac{\pl^2}{\pl t_2\pl s_1}\log \tau_{k_1 ,k_2 }}
      {\frac{\pl^2}{\pl t_1\pl s_1}\log \tau_{k_1 ,k_2 }}
  \label{12}\\
    -\frac{\pl}{\pl s_1}\log \frac{\tau_{k_1 +1,k_2 }}
                           {\tau_{k_1 -1,k_2 }}
                          & =&     \frac
      {\frac{\pl^2}{\pl t_1\pl s_2}\log \tau_{k_1 ,k_2 }}
      {\frac{\pl^2}{\pl t_1\pl s_1}\log \tau_{k _1,k_2 }}
  \label{13}\eea

\end{corollary}

\proof  The bilinear identity (\ref{Psi-bilinear
identity}) yields nine identities, which are all
equivalent, upon relabeling indices, to the tau-function
bilinear identity (\ref{tau-bilinear identity}).
Introducing standard shifts in the residue formulae
  \bean
        t&\longmapsto& t-a ~~~
        t' \longmapsto  t +a\\
        s&\longmapsto& s-b~~~
        s' \longmapsto  s +b\\
        u&\longmapsto& u-c ~~~
        u' \longmapsto  u +c,
        \eean
and using Taylor's theorem, identity (\ref{tau-bilinear
identity}) is equivalent to
\bea & &\sum_{j=0}^{\iy}\gs_{\ell_1+\ell_2-k_1-k_2+j-1}
(-2a)\gs_j(\tilde\pl_t)e^{\sum_1^{\iy}(a_k\frac{\pl}{\pl
t_k}+b_k \frac{\pl}{\pl s_k}+c_k\frac{\pl}{\pl u_k})}
\tau_{\ell_1\ell_2}\circ\tau_{k_1k_2}  \no\\
\no\\
&+&\sum_{j=0}^{\iy}\gs_{k_1-\ell_1+1+j}
(-2b)\gs_j(\tilde\pl_s)e^{\sum_1^{\iy}(a_k\frac{\pl}{\pl t_k}+b_k
\frac{\pl}{\pl s_k}+c_k\frac{\pl}{\pl u_k})}
\tau_{\ell_1-1,\ell_2}\circ\tau_{k_1+1,k_2} \no\\
\no\\
&+&\sum_{j=0}^{\iy}\gs_{k_2-\ell_2+1+j}
(-2c)\gs_j(\tilde\pl_u)e^{\sum_1^{\iy}(a_k\frac{\pl}{\pl
t_k}+b_k \frac{\pl}{\pl s_k}+c_k\frac{\pl}{\pl u_k})}
\tau_{\ell_1,\ell_2-1}\circ\tau_{k_1,k_2+1}=0.\no\\
 \label{5}
\eea
Taylor expanding in $a,b,c$ and setting in equation
(\ref{5}) all $a_i,b_i,c_i=0$, except $a_{j+1}$, and
also setting $\ell_1=k_1+2$, $\ell_2=k_2$, equation
(\ref{5}) becomes

$$
a_{j+1}\left(-2\gs_j(\tilde\pl_t)\tau_{k_1+2,k_2}
\circ\tau_{k_1,k_2}+\frac{\pl^2}{\pl
s_1\pl
t_{j+1}}\tau_{k_1+1,k_2}\circ\tau_{k_1+1,k_2}\right)+
{\bf O}(a^2_{j+1})=0,
$$
and so the coefficient of $a_{j+1}$ must vanish
identically, yielding equation (\ref{bilinear
identity1}), setting $k_1\rg k_1-1$. Setting in equation
(\ref{5}) all $a_i,b_i,c_i=0$, except $b_{j+1}$, and
 $\ell_1=k_1$, $\ell_2=k_2$, the coefficient
of $b_{j+1}$ in equation (\ref{5}) yields equation
(\ref{bilinear identity2}). Specializing equation
(\ref{bilinear identity1}) to $j=0$ and 1 respectively
yields (since $\gs_1(t)=t_1$ implies
$\gs_1(\tilde\pl_t)=\frac{\pl}{\pl t_1}$; also $\gs_0=
1$):
   $$
\frac{\pl^2\log \tau_{k_1 ,k_2 }}{\pl t_1\pl s_1} =
   \frac{\tau_{k_1 +1,k_2 }\tau_{k_1 -1,k_2 }}
         {\tau_{k_1 ,k_2 }^2}
$$
and
$$
\frac{\pl^2}{\pl s_1\pl t_2}\log \tau_{k_1k_2} =
\frac{1}{\tau_{k_1k_2}^2} \left[ \left(\frac{\pl }{ \pl
t_1}
  \tau_{k_1+1,k_2}\right)
   \tau_{k_1-1,k_2}
   -
   \tau_{k_1+1,k_2}
    \left(\frac{\pl }{ \pl t_1}
  \tau_{k_1-1,k_2}\right)
   \right]
  .$$
  Upon dividing the second equation by the first, we find equation
  (\ref{12}) and similarly equation (\ref{13}) follows from equation
  (\ref{bilinear
identity2}). \qed


\section{Virasoro constraints for the integrable
deformations}

%
%
 Given the Heisenberg and Virasoro operators, for $m\geq -1,~ k\geq 0$:
   \bean \BJ^{(1)}_{m,k  }(t)&=&
   \frac{\pl}{\pl t_{m}}+(-m)t_{-m}+k\dt_{0,m}\\
  \BJ^{(2)}_{m,k  }(t)&=&\frac{1}{2}
  \left( \sum_{i+j=m}\frac{\pl^2}{\pl t_i\pl t_j}
   +2\sum_{i\geq 1}it_i\frac{\pl}{\pl t_{i+m}}
     + \sum_{i+j=-m} it_i jt_j \right)
    \\
    &&~~~~~ +(k +\frac{m+1}{2})\left(\frac{\pl}{\pl
     t_m}+(-m)t_{-m}\right)+\frac{k(k+1) }{2}\dt_{m0}
 , \eean
 we now state:

\begin{theorem}
The integral $\tau_{k_1k_2}(t,s,u;\beta;E)$, as defined
in (\ref{1.5})
satisfies

  \be
   {\cal B} _m \tau_{k_1,k_2} =\BV^{k_1,k_2}_m \tau_{k_1,k_2}  ~~\mbox{for}~~
    m\geq -1,
\label{Virasoro constraints}\ee
where ${\cal B}_m$ and $\BV_m$ are differential
operators:
$$
{\cal B}_m=\sum_1^{2r} b_i^{m+1} \frac{\pl}{\pl b_i},
~\mbox{for}~
 E=\bigcup_1^{2r}~[b_{2i-1},b_{2i}] \subset
 \BR
$$
and
  $$
  \BV^{k_1k_2}_m:= \left\{\begin{array}{l}
   \BJ^{(2)}_{m,k_1+k_2}(t)-(m+1)
  \BJ^{(1)}_{m,k_1+k_2}(t)
   \\ \\
  +\BJ^{(2)}_{m,k_1 }(-s)+a
  \BJ^{(1)}_{m+1,k_1}(-s)+(1-2\beta)
  \BJ^{(1)}_{m+2,k_1}(-s)
  \\  \\
   +\BJ^{(2)}_{m,k_2 }(-u)-a
  \BJ^{(1)}_{m+1,k_2}(-u)+(1+2\beta)
  \BJ^{(1)}_{m+2,k_2}(-u)
  \end{array}\right\}
   $$

\end{theorem}
We state the following lemmas:

\begin{lemma} {\em (Adler-van Moerbeke \cite{AvM2})}
 Given
  $$\rho=e^{-V}~~\mbox{with}~~
  -\frac{\rho^{\prime}}{\rho}=V'=\frac{g}{f}
  =\frac{\sum_0^{\iy}\beta_iz^i }{
\sum_0^{\iy}\alpha_iz^i}
,$$
the integrand
 $$
  dI_n(x):= \Dt_n(x) \prod_{k=1}^n
\left(e^{\sum_1^{\iy}t_i x_k^i}\rho(x_k)dx_k\right),$$
  satisfies the following variational
formula:
  \be\left.\frac{d}{d\vr}dI_n (x_i\mapsto
x_i+\vr f(x_i)x_i^{m+1} )\right|_{\vr=0}
=\sum^{\iy}_{\ell=0}
  \left(\al_{\ell}~ \BJ^{(2)}_{m+\ell,n}
-\beta_{\ell}~ \BJ^{(1)}_{m+\ell +1,n}\right)dI_n.
\label{variational1}\ee
The contribution coming from $\prod_{ 1}^{n} dx_j$ is
given by \be
  \sum_{\ell=0}^{\iy} \al_{\ell} (\ell+m+1)
  \BJ^{(1)}_{m+\ell,n} dI_n.
 \label{variational2} \ee

\end{lemma}

\begin{lemma}
Setting
  \bean\lefteqn{ dI_n =\Dt_{n }(x,y)
   \prod_{j=1}^{k_1}e^{\sum_1^{\iy}t_i x_j^i  }
   \prod_{j=1}^{k_2}e^{\sum_1^{\iy}t_i y^i_j  }
   }\\
   & &~~~~~~ \left(\Dt_{k_1}(x)\prod^{k_1}_{j=1}
  e^{-\frac{x_j^2}{2}+ax_j+\beta
x_j^2 }
e^{-\sum_1^{\iy}s_ix^i_j}dx_j\right)\\ &&~~~~~~
\left(\Dt_{k_2}(y)\!  \!   \prod^{k_2}_{j=1}
e^{-\frac{y_j^2}{2}-ay_j-\beta y_j^2}
e^{-\sum_1^{\iy}u_iy^i_j}dy_j\right)
 \eean
   The following variational formula holds for
$m\geq -1$:
  \bea  \left.\frac{d}{d\vr}dI_n \left(
\begin{array}{l}x_i\mapsto x_i+\vr x_i^{m+1}\\
   y_i\mapsto y_i+\vr y_i^{m+1} \end{array}
  \right )\right|_{\vr=0}
  &=&\BV^{k_1,k_2}_m (dI_n ).
 \label{variational3} \eea

\end{lemma}

\proof The variational formula (\ref{variational3}) is
an immediate consequence of applying the variational
formula (\ref{variational1}) separately to the three
factors of $dI_{n}$, and in addition applying formula
(\ref{variational2}) to the first factor, to account for
the fact that
$\prod^{k_1}_{j=1}dx_j\prod^{k_2}_{j=1}dy_j$ is missing
from the first factor.
\qed

\medskip\noindent{\it Proof of Theorem 3.1:\/} Formula (\ref{Virasoro constraints}) follows immediately
from formula
(\ref{variational3})
, by taking into account the variation of $\pl E$ under
the change of coordinates.\qed

%
%


Using the identity, valid when acting on
$\tau_{k_1k_2}(t,s,u;\beta;E)$,
     $$
     \frac{\pl}{\pl t_n}=-\frac{\pl}{\pl s_n}-\frac{\pl}{\pl u_n}
,$$ one obtains by explicit computation for $m\geq -1$,
  $$
  \BV^{k_1k_2}_m:= \left\{\begin{array}{l}
   \BJ^{(2)}_{m,k_1+k_2}(t)-(m+1)
  \BJ^{(1)}_{m,k_1+k_2}(t)
   \\ \\
  +\BJ^{(2)}_{m,k_1 }(-s)+a
  \BJ^{(1)}_{m+1,k_1}(-s)+ (1-2\beta)
  \BJ^{(1)}_{m+2,k_1}(-s)
  \\  \\
   +\BJ^{(2)}_{m,k_2 }(-u)-a
  \BJ^{(1)}_{m+1,k_2}(-u) +(1+2\beta)
  \BJ^{(1)}_{m+2,k_2}(-u)
  \end{array}\right\}
    $$

   \bean
   &=&\frac{1}{2}
   \sum_{i+j=m}
   \left(
    \frac{\pl^2}{\pl t_i\pl t_j}
    +\frac{\pl^2}{\pl s_i\pl s_j}
    +\frac{\pl^2}{\pl u_i\pl u_j}\right)
+
   \sum_{i\geq 1}
    \left(it_i\frac{\pl}{\pl t_{i+m}}
         +is_i\frac{\pl}{\pl s_{i+m}}
         +iu_i\frac{\pl}{\pl u_{i+m}}
    \right)
    \\
    &&
    +(k_1+k_2)\left(\frac{\pl}{\pl
     t_m}+(-m)t_{-m}\right)
     -k_1 \left(\frac{\pl}{\pl
     s_m}+(-m)s_{-m}\right)
     -k_2 \left(\frac{\pl}{\pl
     u_m}+(-m)u_{-m}\right)
   \\
   &&
   +(k_1^2+k_1k_2+k_2^2)\dt_{m0}+a(k_1-k_2)\dt_{m+1,0}
   +\frac{m(m+1)}{2}(-t_{-m}+s_{-m}+u_{-m})
  \\
  && - \frac{\pl}{\pl
     t_{m+2}}
     +a\left(-\frac{\pl}{\pl
     s_{m+1}}+ \frac{\pl}{\pl
     u_{m+1}}+(m+1)(s_{-m+1}-u_{-m+1})\right)
     \\
     &&
     + 2\beta \left(\frac{\pl}{\pl
     u_{m+2}}-\frac{\pl}{\pl
     s_{m+2}} \right)\eean

The following identities, valid when acting on
$\tau_{k_1k_2}(t,s,u;\beta;E)$, will also be used:
   \bean
\begin{array}{llll}
\frac{\pl}{\pl s_1}=-\frac{1}{2}\left(\frac{\pl}{\pl
t_1}+\frac{\pl}{\pl a}\right)& & &\frac{\pl}{\pl
s_2}=-\frac{1}{2}\left(\frac{\pl}{\pl t_2}
  +\frac{\pl}{\pl \beta}\right)\\
\\
\frac{\pl}{\pl u_1}=-\frac{1}{2}\left(\frac{\pl}{\pl
t_1}-\frac{\pl}{\pl a}\right)& & &\frac{\pl}{\pl
u_2}=-\frac{1}{2}\left(\frac{\pl}{\pl
t_2}-\frac{\pl}{\pl \beta}\right)
\end{array}
\eean


\begin{corollary} 
The $tau$-function $\tau=\tau_{k_1,k_2}(t,s,u;\beta;E)$
satisfies
     the following differential identities, with
     ${\cal B}_m=\sum_1^{2r} b_i^{m+1} \frac{\pl}{\pl b_i}$:

\bean
-B_{-1}\tau&=&\left(\frac{\pl}{\pl t_1}
 -2\beta \frac{\pl }{\pl a}
\right)\tau
 \\
 & &-\sum_{i\geq 2}\left(it_i\frac{\pl}{\pl
t_{i-1}} +is_i\frac{\pl}{\pl s_{i-1}}+iu_i\frac{\pl}{\pl
u_{i-1}}\right)\tau\\
& &
+a(k_2-k_1)\tau+(k_1s_1+k_2u_1-(k_1+k_2)t_1)\tau\\ \\
\frac{1}{2}\left(B_{-1}-\frac{\pl}{\pl a}\right)\tau&=&
 \left(\frac{\pl }{\pl s_{1}} +\beta
\frac{\pl }{\pl a} \right)\tau\\
 & &+\frac{1}{2}\sum_{i\geq 2}
\left(it_i\frac{\pl}{\pl t_{i-1}} +is_i\frac{\pl}{\pl
s_{i-1}}+iu_i\frac{\pl}{\pl
u_{i-1}}\right)\tau\\
& &
+\frac{a}{2}(k_1-k_2)\tau
 +\frac{1}{2}((k_1+k_2)t_1-k_1s_1-k_2u_1)\tau
  \\ \\
-\left(B_0-a\frac{\pl}{\pl a}\right)\tau&=&
\frac{\pl\tau}{\pl
t_2}-(k_1^2+k_2^2+k_1k_2)\tau\\
& &-2\beta \frac{\pl\tau}{\pl \beta}-\sum_{i\geq
1}\left(it_i\frac{\pl}{\pl t_i}+is_i\frac{\pl}{\pl
s_i}+iu_i\frac{\pl}{\pl u_i}\right)\tau\\
\frac{1}{2}\left(B_0-a\frac{\pl}{\pl a}-\frac{\pl}{\pl
\beta}\right)\tau&=&\frac{\pl\tau}{\pl
s_2}+\frac{1}{2}(k_1^2+k_2^2+k_1k_2)\tau\\
& &+\beta\frac{\pl\tau}{\pl
\beta}+\frac{1}{2}\sum_{i\geq 1}\left(it_i\frac{\pl}{\pl
t_i}+is_i\frac{\pl}{\pl s_i}+iu_i\frac{\pl}{\pl
u_i}\right)\tau \eean \be \label{8}\ee
\end{corollary}

\begin{corollary} On the locus $\LR=\{t=s=u=0,
\beta=0\}$, the function $f=\log
\tau_{k_1k_2}(t,s,u;\beta;E)$ satisfies the following
differential identities:

  \bean
\frac{\pl f}{\pl t_1} &=&-\BB_{-1}f
    +a(k_1-k_2)  
\no\\
\frac{\pl f}{\pl s_1}
&=&\frac{1}{2}\left(\BB_{-1}-\frac{\pl  }{\pl
a}\right)f+\frac{a}{2}(k_2-k_1)  \no\\
  \frac{\pl f}{\pl
t_2} &=&\left(-\BB_0+a\frac{\pl}{\pl
a}\right)f+k_1^2+k_1k_2+k_2^2 
 \eean\bea \frac{\pl f}{\pl
s_2} &=&\frac{1}{2}\left(\BB_{0}-a\frac{\pl }{\pl
a}-\frac{\pl}{\pl
\beta}\right)f-\frac{1}{2}(k_1^2+k^2_2+k_1k_2)
\label{16}\\  \no\\
 2\frac{\pl^2f}{\pl t_1\pl s_1}
&=&\BB_{-1}\left(\frac{\pl}{\pl a}-
\BB_{-1} \right)f-2k_1\no\\
 \no\\
2\frac{\pl^2f}{\pl t_1\pl s_2} &=&\left( a
\frac{\pl}{\pl a} +\frac{\pl}{\pl \beta} -\BB_0
+1\right)\BB_{-1}f-2\frac{\pl f}{\pl a}
  -2a(k_1-k_2)\no\\
\no\\
2\frac{\pl^2f}{\pl t_2\pl s_1} &=& \frac{\pl}{\pl
a}(\BB_0-a\frac{\pl}{\pl a}+a\BB_{-1})f
-\BB_{-1}(\BB_0\!-\!1)f-2a(k_1\!-\!k_2)
\label{20}
    \eea


\end{corollary}

    \proof Upon dividing equations (\ref{8}) by $\tau$ and restricting
    to the locus $\LR$, equations (\ref{16}) follow immediately. The
    essence of deriving (\ref{20}) is that the
     Virasoro operators $\BV_n$
    and the boundary operators $\BB_m$ commute.
     To derive, say, the first equation in the list
    (\ref{20}), rewrite the two first equations
    of (\ref{8}) as
    \bean
    -\BB_{-1}f&=&\frac{\pl f}{\pl t_1} +a(k_2-k_1)+L_1(f)+\ell_1\\
    \frac{1}{2}\left(\BB_{-1}-\frac{\pl}{\pl a}\right)f&=&\frac{\pl f}{\pl s_1}
    +\frac{1}{2}{a(k_1-k_2)} +L_2(f)+\ell_2
    \eean
    where $L_i$ are linear operators vanishing on $\LR$ and the
    $\ell_i$ are functions vanishing on $\LR$. This
    yields:
\bean \lefteqn{ (-{\cal B}_{-1})\frac{1}{2}\left({\cal
B}_{-1}-\frac{\pl}{\pl
a}\right)f\Big\vert_{\LR}}\\
&=&\left(\frac{\pl}{\pl s_1}+\beta\frac{\pl}{\pl
a}\right)(-{\cal
B}_{-1}f)\Big\vert_{\LR}\\
\\
& &\qquad\qquad +\frac{1}{2}\sum_{i\geq
2}\left(it_i\frac{\pl}{\pl t_{i-1}}+is_i\frac{\pl}{\pl
s_{i-1}}+iu_i \frac{\pl}{\pl u_{i-1}}\right)(-{\cal
B}_{-1}f)\Big\vert_{\LR}\\
\\
&=&\frac{\pl}{\pl s_1}(-{\cal B}_{-1}f)\Big\vert_{\LR}\\
\\
&=&\frac{\pl}{\pl s_1}\left.\left(\begin{array}{l}
\left(\frac{\pl}{\pl t_1}-2\beta\frac{\pl}{\pl a}\right)f+a(k_2-k_1)\\
-\displaystyle{\sum_{i\geq 2}\left(it_i\frac{\pl}{\pl
t_{i-1}}+is_i\frac{\pl}{\pl s_{i-1}}+
iu_i\frac{\pl}{\pl u_{i-1}}\right)f}\\
+~(k_1s_1+k_2u_1-(k_1+k_2)t_1)
\end{array}\right)\right\vert_{\LR}=\frac{\pl^2}{\pl s_1\pl t_1}f+k_1
.\eean
  The other
identities (\ref{20}) can be obtained in a similar
way.\qed


\section{A PDE for the Gaussian ensemble with external
source}  

   \medskip\noindent{\it Proof of Theorem 0.1:\/}
   First observe that, with $n=k_1+k_2$,
   $$ \BP_{n}(a; E)=\frac{1}{Z_{n}}
  \int_{{\cal H}_{n}(E)}
  e^{-\Tr(\frac{1}{2} M^2-AM) }dM
   =\frac{\tau_{k_1k_2}(t,s,u;\beta;E)}
         {\tau_{k_1k_2}(t,s,u;\beta;\BR)}
   \Bigr|_{t=s=u=\beta=0}$$
An explicit computation over the whole range
yields:
   \bean
 \lefteqn{  \tau_{k_1k_2}(t,s,u;\beta;\BR)
   \Bigr|_{t=s=u=\beta=0}}\\
   &=&
    \frac{1}{k_1!k_2!}\int_{\BR^{n}}\Dt_{n}(x,y)
 \left(\Dt_{k_1}(x)\prod^{k_1}_{i=1}e^{-\frac{x_i^2}{2}+a
x_i}dx_i\right)\\&&~~~~~~~~~~~~~~~~~~~~~~~~~~
\left(\Dt_{k_2}(y)\prod^{k_2}_{i=1}
e^{-\frac{y^2_i}{2}-ay_i}dy_i\right)\\
&=&
c_{k_1k_2}a^{k_1k_2}e^{(k_1+k_2)a^2/2}. \eean
 This is obtained from the representation (\ref{1.5})
  in terms of moments, which themselves
   are Gaussian integrals,as shown in Appendix 1. From this formula, it follows
   that
  $$
\log \tau_{k_1k_2}(t,s,u;\beta;\BR)
   \Bigr|_{t=s=u=\beta=0}=\frac{k_1+k_2}{2}a^2+k_1k_2\log a+C_{k_1k_2},
$$
where $c_{k_1k_2}$ and $C_{k_1k_2}$ are constants
depending on ${k_1,k_2}$ only. It follows that
  \be
  \log \BP_{n}(a;E)= \log \tau_{k_1k_2}(0,0,0;0;E)
   -\frac{k_1+k_2}{2}a^2-k_1k_2\log a-C_{k_1k_2}
   \label{4.1}\ee
Thus we need to concentrate on
$\tau_{k_1k_2}(t,s,u;\beta;E)$, which, by Theorem 2.1,
satisfies the bilinear identity (\ref{Psi-bilinear
identity}) and thus the identities (\ref{12}) and
(\ref{13}) of Corollary 2.2:
%
%
\bea   \frac{\pl}{\pl t_1}\log \frac{\tau_{k_1  +1,k_2
}}
                           {\tau_{k_1  -1,k_2  }}
                           &=&\frac
      {\frac{\pl^2}{\pl t_2\pl s_1}\log \tau_{k _1  k_2  }}
      {\frac{\pl^2}{\pl t_1\pl s_1}\log \tau_{k_1   k_2  }}
  \no\\
    \frac{\pl}{\pl s_1}\log \frac{\tau_{k_1  +1,k_2  }}
                           {\tau_{k_1  -1,k_2 }}
                         &=& -\frac
      {\frac{\pl^2}{\pl t_1\pl s_2}\log \tau_{k_1   k_2  }}
      {\frac{\pl^2}{\pl t_1\pl s_1}\log \tau_{k_1   k_2  }},
   \label{a}\eea
whereas the first two Virasoro equations (\ref{16})
yield, specializing to the locus $\LR=\{t=s=u=0,
\beta=0\}$ and the indices $k_1\pm 1, k_2$,
 \bea
 \frac{\pl}{\pl t_1}\log \frac{\tau_{k_1  +1,k_2  }}
                           {\tau_{k_1  -1,k_2  }}
                           &=&
       -\BB_{-1}\log \frac{\tau_{k_1  +1,k_2  }}
                           {\tau_{k_1  -1,k_2  }}+2a
                          \no \\
 \frac{\pl}{\pl s_1}\log \frac{\tau_{k_1 +1,k_2  }}
                           {\tau_{k_1  -1,k_2  }}
                          & =&
       \frac{1}{2}\left(\BB_{-1}
        -\frac{\pl  }{\pl a}\right)
       \log\frac{\tau_{k_1  +1,k_2  }}
                           {\tau_{k_1 -1,k_2  }}- a
  \label{b}
  \eea
From these three equations, the expression
$\log\frac{\tau_{k_1 +1,k_2 }} {\tau_{k_1 -1,k_2  }}$
can be eliminated, by first subtracting the first
equations in (\ref{a}) and (\ref{b}) and then the second
equations in (\ref{a}) and (\ref{b}). Subsequently one
eliminates $\log\frac{\tau_{k_1 +1,k_2 }} {\tau_{k_1
-1,k_2 }}$ from the equations thus obtained, yielding
$$
   \frac{1}{2}\left(\BB_{-1}- \frac{\pl}{\pl a} \right)
  \left( \frac
      {\frac{\pl^2}{\pl t_2\pl s_1}\log \tau_{k_1k_2}}
      {\frac{\pl^2}{\pl t_1\pl s_1}\log \tau_{k_1k_2}}
  -2a\right)=
    \BB_{-1}
   \left( \frac
      {\frac{\pl^2}{\pl t_1\pl s_2}\log \tau_{k_1k_2}}
      {\frac{\pl^2}{\pl t_1\pl s_1}\log \tau_{k_1k_2}}
-a\right)$$
  or equivalently
\be \BB_{-1}\!\left( \frac{\left(\frac{\pl^2  }{\pl
t_2\pl s_1}
               \!\! -\!\!2\frac{\pl^2  }{\pl t_1\pl s_2}\right)
                 \log\tau_{k_1k_2  }}
                {\frac{\pl^2  }{\pl t_1\pl s_1}
                \log\tau_{k_1 k_2}}
\right)\!
    -\frac{\pl}{\pl a}\!\!\left(
    \frac{\left(\frac{\pl^2  }{\pl t_2\pl s_1}
                \!\!-\!\!  2a\frac{\pl^2  }{\pl t_1\pl s_1}\right)
                \log\tau_{k_1 k_2}}
                {\frac{\pl^2  }{\pl t_1\pl s_1}
                \log\tau_{k_1 k_2 }}
\right)=0 . 
\label{4.4}\ee
Using the Virasoro relations (\ref{20}), one obtains
along the locus $\LR=\{t=s=u=0, \beta=0\}$:
 \bea
  4\frac{\pl^2  }{\pl t_1\pl s_1}\log\tau_{k_1k_2}&=:&
   F^+
  \no\\
2\left(\frac{\pl^2  }{\pl t_2\pl s_1}
  -2 \frac{\pl^2  }{\pl t_1\pl s_2}\right)
   \log\tau_{k_1k_2}
  &=&  H_1^+-2\BB_{-1}\frac{\pl  }{\pl \beta}\log\tau_{k_1k_2}
  \no\\
2\left(\frac{\pl^2  }{\pl t_2\pl s_1}
  -2a \frac{\pl^2  }{\pl t_1\pl s_1}\right)\log\tau_{kk}
  &=&  H_2^+
  \label{4.5}\eea
where, using the identity (\ref{4.1}), along the locus
$\LR=\{t=s=u=0, \beta=0\}$,
 {\footnotesize \bean
       F^+&:=&  2\BB_{-1}(\frac{\pl}{\pl
   a}\!-\!\BB_{-1})\log \tau_{k_1k_2}-4k_1
   =2\BB_{-1}(\frac{\pl}{\pl
   a}\!-\!\BB_{-1})\log \BP_{n}-4k_1 
   \\
        H_1^+&:=& \frac{\pl}{\pl a}
      \left( \BB_{0} -a\frac{\pl}{\pl a}-a \BB_{-1}
      \right)\log \tau_{k_1k_2}
       +\left(\BB_{0}\BB_{-1}+4\frac{\pl}{\pl a}\right)
       \log \tau_{k_1k_2}+2a(k_1-k_2)
      \\
      &=& \frac{\pl}{\pl a}
      \left( \BB_{0} -a\frac{\pl}{\pl a}-a \BB_{-1}
      \right)\log \BP_{n}
       +\left(\BB_{0}\BB_{-1}+4\frac{\pl}{\pl a}\right)
       \log \BP_{n}
      +4a k_1 +\frac{4k_1k_2}{a}\\
      \\ \hspace{-1cm}
       H^+_2
      &:=&\!\!\!
      \frac{\pl}{\pl a}
      \left( \BB_{0} -a\frac{\pl}{\pl a}-a \BB_{-1}
      \right)\log \tau_{k_1k_2}
        +\left(2a\BB_{-1}
       \!-\!\BB_{0}
     \! +\!2\right)\BB_{-1}\log \tau_{k_1k_2} \!+2a(k_1\!+\!k_2)\\
      &=&\!\!\!\frac{\pl}{\pl a}
      \left( \BB_{0} -a\frac{\pl}{\pl a}-a \BB_{-1}
      \right)\log \BP_{n}
        +\left(2a\BB_{-1}
       -\BB_{0}
      +2\right)\BB_{-1}\log \BP_{n},
   \eean}
confirming (\ref{0.7}).
Notice that the expressions above do not contain
partials in $\beta$, except for the $\beta$-partial
appearing in the second expression of (\ref{4.5}).
Putting these expressions (\ref{4.5}) into (\ref{4.4})
yields
   \bea
   \left\{\left.\BB_{-1}\frac{\pl }{\pl \beta}
    \log \tau_{k_1k_2}\right|_{\LR}~,
     F^+
    \right\}_{\BB_{-1}}
&=&  \left\{
       H_1^+,
       \frac{1}{2}  F^+
       \right\}_{\BB_{-1}}
      -\left\{
       H^+_2,\frac{1}{2}  F^+
     \right\}_{\pl/\pl a}\no\\
     &=:& G^+
  \label{4.6}   \eea
and by involution $a \mapsto -a,~  \beta \mapsto
-\beta,~k_1 \longleftrightarrow k_2$:
      \bea
  -\left\{\left.\BB_{-1}\frac{\pl  }{\pl \beta}
  \log \tau_{k_1k_2}\right|_{\LR}~,   F^-
    \right\}_{\BB_{-1}}
&=&   \left\{
       H_1^-,
       \frac{1}{2}  F^-
       \right\}_{\BB_{-1}}
      -\left\{
       H^-_2,\frac{1}{2}  F^-
     \right\}_{-\pl/\pl a}\no\\
     &=:& G^-
   \label{4.7}   \eea
   where
  $$
  F^-=  F ^+\Bigr|_{\begin{array}{l}
              a\rightarrow -a \\
               k_1 \leftrightarrow k_2
              \end{array}}
              ,
   H_i^-=  H_i ^+\Bigr|_{\begin{array}{l}
              a\rightarrow -a \\
              k_1 \leftrightarrow k_2
              \end{array}}.
  $$
  Remember the change of variables $a \mapsto -a,~  \beta \mapsto
-\beta,~k_1 \longleftrightarrow k_2$ acts on the
operators,
              since $\tau_{k_1k_2}$ is invariant under this
              change; see (\ref{1.6}).

\noindent Equations (\ref{4.6}) and (\ref{4.7}) yield a
linear
  system of equations in
  $$
   \BB_{-1}\frac{\pl \log \tau_{k_1k_2} }{\pl \beta} ~~~\mbox{and}
   ~~~\BB_{-1}^2\frac{\pl \log \tau_{k_1k_2} }{\pl \beta}
  $$
from which
  \bean
   \BB_{-1}\frac{\pl \log \tau_{k_1k_2}}{\pl \beta}&=&
    \frac{G^- F^+ + G^+F^-}
    {-F^-(\BB_{-1}F^+)+F^+(\BB_{-1}F^-)}
\\
 \BB_{-1}^2\frac{\pl \log \tau_{k_1k_2}}{\pl \beta}&=&
    \frac{G^- (\BB_{-1}F^+) + G^+(\BB_{-1}F^-)}
    {-F^-(\BB_{-1}F^+)+F^+(\BB_{-1}F^-)}
\eean
  Subtracting the second equation from
$\BB_{-1}$ of the first equation yields the following:
  \bean&&
\Bigl(F^+\BB_{-1}G^-+F^-\BB_{-1}G^+ \Bigr)
\Bigl(F^+\BB_{-1}F^- -F^-\BB_{-1}F^+ \Bigr) \\
  &&-
  \Bigl(F^+ G^- +F^- G^+
\Bigr) \Bigl(F^+\BB_{-1}^2F^- -F^-\BB_{-1}^2F^+ \Bigr)
=0,
 \eean
 establishing Theorem 0.1.\qed


\section{A PDE for the Pearcey transition probability}

From the Karlin-McGregor formula for non-intersecting
Brownian motions $x_j(t)$, we have:

\bean
 \lefteqn{\BP\left( \mbox{all}~ x_i(t)
  \in E,~1\leq i\leq n \left|
 \begin{array}{l}
  \mbox{given $x_i(0)=\gamma_i$}\\
  \mbox{given $x_i(1)=\delta_i$}  \end{array}
  \right.\right)}\\
  &=&
  \int_{E^n}\frac{1}{Z_{n}}\det(p(t;\gamma_i,x_j))_{1\leq i,j\leq
n}\det(p(1-t;x_{i'},\delta_{j'}))_{1\leq i',j'\leq n}
 \prod_1^n dx_i\eean
  for the Brownian motion kernel
  $$
  p(t,x,y):= \frac{1}{\sqrt{ \pi t}
  }~e^{-\frac{(y-x)^2}{ t}}.
  $$
Aptekarev, Bleher and Kuijlaars introduce in
\cite{AptBleKui} a change of variables transforming the
Brownian motion problem into the Gaussian random
ensemble with external source. For $E:=
\bigcup^r_{i=1}[b_{2i-1},b_{2i}]$, we have, using this
change of variables,
$$
x _i=x'_i\sqrt{\frac{t(1-t)}{2}} ~~~\mbox{and}~~~ y
_i=y'_i\sqrt{\frac{t(1-t)}{2}},
$$
in equality $\stackrel{*}{=}$,
  \bean
\lefteqn{ \BP_0^{\pm a}\left(\mbox{all $x_j(t)\in
E$}\right)
  }
   \\
&:=& \BP \left(\begin{tabular}{c|c}
   & all $x_j(0) =0$\\
  all $x_j(t)\in E$  & $k$ left paths end up at $-a$ at time $t=1$,\\
   &
$k$ right paths end up at $+a$ at time $t=1$
\end{tabular}\right)
\eean\bean   &=&\lim_{\begin{array}{c}
  \mbox{\small all}~\gamma_i\rightarrow 0
  \\
  \delta_1, \ldots , \delta_k\rightarrow -a\\
  \delta_{k+1}, \ldots , \delta_{2k}\rightarrow a
  \end{array}}\\&&~~~
  \int_{E^n}\frac{1}{Z_{n}}\det(p(t;\gamma_i,x_j))_{1\leq i,j\leq
n}\det(p(1-t;x_{i'},\delta_{j'}))_{1\leq i',j'\leq n}
  \prod_1^n dx_i\\
  &=& \frac{1}{Z_{n}}
  \int_{E^{n}} \Dt_{n}(x ,y )
\left(\Dt_k(x )\prod^k_{i=1}e^{ -\frac{x_i^{
2}}{t(1-t)}+\frac{2ax _i}{1-t}}dx _i\right)
%
\!\! \left(\Dt_k(y )\prod^k_{i=1}e^{ -\frac{y_i^{
 2}}{t(1-t)}-\frac{2ay _i}{1-t}}dy _i\right)
 \\
 &\stackrel{*}{=}&
  \frac{1}{Z'_{n}}
  \int_{\left(E\sqrt{\frac{2}{t(1-t}}\right)^{n}} \Dt_{n}(x' ,y' )
\left(\Dt_k(x' )\prod^k_{i=1}e^{ -\frac{x_i^{\prime
2}}{2}+ { a\sqrt{\frac{2t }{1-t}}}x' _i }dx' _i\right)
\\
&& ~~~~~~~~~~~~~~~~~~~~~~~~~~~~~~~~~~~~~ \left(\Dt_k(y'
)\prod^k_{i=1}e^{ -\frac{y_i^{\prime
 2}}{2}- a\sqrt{\frac{2t }{1-t}}y' _i }dy' _i\right)
 \\
  &=&\BP_{n}\left(a \sqrt{\frac{2t }{1-t}};
b_1\sqrt{\frac{2}{t(1-t)}},\ldots, b_{2r
}\sqrt{\frac{2}{t(1-t)}}\right)
\eean with $\BP_n$ of Theorem 0.1, using (\ref{1.1}) and
(\ref{1.3}), with $k=k_1=k_2$.
%
%
%
Setting \be
 e^{g(t)}:=\sqrt{\frac{2t}{1-t}}\mbox{ and }
 e^{h(t)}:=\sqrt{\frac{2}{t(1-t)}}
  \label{5.1}\ee

\be \tilde\BB_k=\sum_1^{2r}v_i^{k+1}\frac{\pl}{\pl
v_i},\qquad\BB_k=\sum_1^{2r}b_i^{k+1}\frac{\pl}{\pl
b_i},\label{5.2}\ee we find

 \be
\BP_0^{\pm
a}(t;b_1,\ldots,b_{2r})=\BP_{n}(ae^{g(t)};b_1e^{h(t)},\ldots,
 b_{2r}e^{h(t)})=
\BP_{n}(u;v_1,\ldots,v_{2r})\Big\vert_{{u=ae^{g(t)}}
 \atop{v=be^{h(t)}}}.
\label{5.3}
 \ee
 From Theorem 0.1, it follows that
 $\BP_{n}(u;v_1,\ldots,v_{2r})$
  satisfies the non-linear equation (\ref{0.5}), with
  $a$ and all $b_i$'s replaced by $u$ and $v_i$
  respectively. In order to find the equation for $\BP_0^{\pm
a}(t;b_1,\ldots,b_{2r})$, one needs to compute the
partial derivatives in $t_i$ and $b_i$ in terms of
partials in $u$ and $v_i$, appearing in equation
(\ref{0.5}) and use the relationship (\ref{5.3}). To be
precise, compute
\be
   \left(\frac{\pl}{\pl
t}\right)^i(\BB_{0})^j(\BB_{-1})^{\ell} ~\BP_0^{\pm
a}~\mbox{with}~ i+j+{\ell}\leq 4 \mbox{~and~}
  i,j,{\ell} \geq 0,
  \label{5.4} \ee
 yielding a system of 34 linear equations in 34 unknowns
   \be
\left(\frac{\pl}{\pl u}\right)^i(\tilde
\BB_{0})^j(\tilde\BB_{-1})^{\ell} ~\BP_{n}~\mbox{with}~
i+j+{\ell}\leq 4 \mbox{~and~}
  i,j,{\ell} \geq 0,
  \label{5.5} \ee
   which one solves.
Notice, one always writes $(\tilde
\BB_{0})^j(\tilde\BB_{-1})^{\ell} $ in that order, using
the commutation relation
$[\tilde\BB_{-1},\tilde\BB_0]=\tilde\BB_{-1}$.
For instance,
\bean
 (\BB_{-1})^j\BP_0^{\pm a}&=&e^{jh(t)}(\tilde\BB_{-1})^j
 \BP_{n},~~~
  (\BB_{0})^j\BP_0^{\pm a}= (\tilde\BB_{0})^j
 \BP_{n}\qquad j=1,\ldots, 4.\\
\\
 \frac{\pl \BP_0^{\pm a}}{\pl t}&=&
  \left(g'(t)u\frac{\pl}{\pl
u}+h'(t)\tilde\BB_0\right)\BP_{n}\\
\\
  \frac{\pl^2 \BP_0^{\pm
a}}{\pl t^2}  &=&\left(g'(t)u\frac{\pl}{\pl
u}+h'(t)\tilde\BB_0
 \right)
\left(g'(t)u\frac{\pl}{\pl u}+h'(t)\tilde\BB_0\right)
 \BP_{n}\\
\\
&&~~~ ~+ \left(g''(t)u\frac{\pl}{\pl
u}+h''(t)\tilde\BB_0\right)\BP_{n} \\
 & &\,\,\,\vdots \eean
 The partials (\ref{5.5}) thus obtained are now
  being substituted into the 4th order
 equation (\ref{0.5}), with $a$ and $b_i$ replaced by $u$ and $v_i$,
 and thus the $\BB_j$ by $\tilde \BB_j$, yielding a new 4th order
 equation involving the partials (\ref{5.4}).

Let now the number of particles $n$ go to infinity,
together with the corresponding scaling (see
\cite{AptBleKui,TW})

  \be
n=2k=\frac{2}{z^4},\quad \pm a=\pm\frac{ 1}{z^2},\quad
b_i=x_iz,\quad t=\frac{1}{2}+sz^2,~~~\mbox{~for~}z\rg 0.
\label{5.6}  \ee
 It is convenient to replace the $\pm$ in
(\ref{5.6}) by the variable $\vr$, which one keeps in
the computation as a variable.
The scaling combined with the change of variables
(\ref{5.3}) leads to the following expressions $u$ and
$v_i$ in terms of $z$:
 \bea
  2k&=& \frac{2}{z^4}\no\\
  u&=&ae^{g(t)}=a\sqrt{\frac{2t}{1-t}}
  = \frac{\vr\sqrt{2}}{z^2}
                    \sqrt{\frac{\frac{1}{2}+sz^2}
                    {\frac{1}{2}-sz^2}}\no\\
  v_i&=&b_ie^{h(t)}=
    b_i\sqrt{\frac{2}{t(1-t)}}=
     x_i\frac{z\sqrt{2}}{\sqrt{\frac{1}{4}-s^2z^4}}
\label{5.7}  \eea So, the question now is to estimate:
{\footnotesize  \bea
\left.\left\{\begin{array}{l}
\Bigl(F^+\tilde\BB_{-1}G^-+F^-\tilde\BB_{-1}G^+ \Bigr)
\!  \Bigl(F^+\tilde\BB_{-1}F^- -F^-\tilde\BB_{-1}F^+ \Bigr) \\
 -
 \Bigl(F^+ G^- +F^- G^+
\Bigr) \Bigl(F^+\tilde\BB_{-1}^2F^-
-F^-\tilde\BB_{-1}^2F^+ \Bigr)
\end{array}\right\}
 \right|_{\begin{array}{l}
                    u\mapsto \frac{\vr\sqrt{2}}{z^2}
                    \sqrt{\frac{\frac{1}{2}+sz^2}
                    {\frac{1}{2}-sz^2}} \\
                    v_i\mapsto x_i\frac{z\sqrt{2}}{\sqrt{\frac{1}{4}-s^2z^4}}
          \\
          n \mapsto  \frac{2}{z^4}\end{array} }
  \label{5.8}\eea
  }
For this, we need to compute the expressions
$F^{\pm},\tilde\BB_{-1} F^{\pm},
\tilde\BB_{-1}^2F^{\pm}, G^{\pm}$ and $\tilde\BB_{-1}
G^{\pm}$ appearing in (\ref{5.8}) in terms of
 { \footnotesize \bea
  \lefteqn{Q_z(s;x_1,\ldots,x_{2r})}
  \no\\
  &:=& \log\BP_{2/z^4}\left(  \frac{\vr\sqrt{2}}{z^2}
                    \sqrt{\frac{\frac{1}{2}+sz^2}
                    {\frac{1}{2}-sz^2}}~;~
                    x_1\frac{z\sqrt{2}}
                    {\sqrt{\frac{1}{4}-s^2z^4}},
                    \ldots,
                    x_{2r}\frac{z\sqrt{2}}
                    {\sqrt{\frac{1}{4}-s^2z^4}}
                     \right)
                     \no\\
     &=& Q(s;x_1,\ldots, x_{2r})+O(z)                ,
  \label{5.10}\eea}
  with
  \be
   Q(s;x_1,\ldots, x_{2r})=\log \det
  \left(I-K_s\chi_{_{E^c}}\right)
   ,\label{5.11}\ee
   as shown in \cite{TW}.
    Without taking a limit on
   $Q_z(s;x_1,\ldots,x_{2r})$ yet, one computes
%
{\footnotesize \bea
  F^{\vr}& =& -\frac{4}{z^4}-\frac{1}{4z^2}\BB_{-1}^2 Q_z+
 \frac{\vr}{4  z}\BB_{-1}\frac{\pl  Q_z}{ \pl
 s }+O(z)
 \no\\ \no \\
  \frac{1}{\sqrt{2}}\tilde\BB_{-1}F^{\vr}
 &=&
 -\frac{1}{16z^3}\BB_{-1}^3 Q_z
 +
\frac{\vr}{16  z^2}\BB_{-1}^2\frac{\pl Q_z}{\pl s}
 -\frac{\vr s}{8   }\BB_{-1}^2\frac{\pl Q_z}{\pl s}
   +  O(z)
\no \\  \no\\
   \tilde\BB_{-1}^2F^{\vr}
 &=&
 -\frac{1}{32z^4}\BB_{-1}^4 Q_z
 +
\frac{\vr}{32  z^3}\BB_{-1}^3\frac{\pl Q_z}{\pl s}
 -\frac{\vr s}{16  z }\BB_{-1}^3\frac{\pl Q_z}{\pl s}
   +  O(1)
  \no\\ \no  \\
  G^{\vr} &=&\frac{3\vr}{8  z^9}\BB_{-1}^3 Q_z
   +\frac{  \vr s}{4  z^7}
    \BB_{-1}^3 Q_z
 \no\\
 &&
-\frac{1}{128z^6} \left[
   (\BB_{-1}\frac{\pl Q_z}{\pl s})(\BB_{-1}^3
   Q_z)+32\BB_0\BB^2_{-1}Q_z
    -  \left(\BB_{-1}^2Q_z+64s\right)\BB_{-1}^2\frac{\pl
    Q_z}{\pl s}
   \right.
   \no\\
   &&
  ~~~~~~~~~~~~~~\left. -64 \BB_{-1}^2Q_z
    +16\frac{\pl^3Q_z}{\pl s^3}
   \right]
   +O(\frac{1}{z^5})
   \no\\ \no\\
\frac{1}{\sqrt{2}} \tilde\BB_{-1} G^{\vr}
 &=&\frac{3\vr}{32
z^{10}}\BB_{-1}^4 Q_z
   +\frac{ \vr s}{16  z^8}\BB_{-1}^4Q_z
 \no\\
 &&
+\frac{1}{512z^7} \left[
 -(\BB_{-1}\frac{\pl Q_z}{\pl s})(\BB_{-1}^4Q_z)
 -32\BB_0\BB_{-1}^3Q_z
   + \left(\BB_{-1}^2Q+64s\right)\BB_{-1}^3\frac{\pl Q_z}{\pl s}
   \right.\no\\
   &&~~~~~~~~~~~~~\left.  +32 \BB_{-1}^3 Q_z
    -16 \BB_{-1}\frac{\pl^3Q_z}{\pl s^3}
   \right]
   +O(\frac{1}{z^6}).\label{5.10}
   \eea }
The formulae needed to obtain the expansions above for $
G^{\vr}$ and $\tilde\BB_{-1}G^{\vr}$ are given in
Appendix 2. From the expressions above one readily
deduces
 {\footnotesize
 \bean
   F^+\tilde\BB_{-1}G^-+F^-\tilde\BB_{-1}G^+
\!\!&=&\!\!- \frac{\sqrt{2}}{64z^{11}}\!
 \left( \!\! \begin{array}{l}
  2(\BB_{-1}\frac{\pl Q_z}{\pl s})(\BB_{-1}^4Q_z )
    \\ \\
    -32(\BB_0-2s\frac{\pl}{\pl s}-1)\BB_{-1}^3Q_z
     \\ \\
     +
      (\BB_{-1}^2\frac{\pl Q_z}{\pl s})(\BB_{-1}^3Q_z )
      -16 \BB_{-1}\frac{\pl ^3Q_z}{\pl s^3}
     \end{array}
      \right)
     \! +\!O(\frac{1}{z^{9}})
 \\
  F^+\tilde\BB_{-1}F^- -F^-\tilde\BB_{-1}F^+
  \!\!&=&\! \vr\frac{\frac{\pl}{\pl s}\BB_{-1}^2
Q_z}{ \sqrt{2}z^{6}}+O(\frac{1}{z^{4}})
\\
   F^+ G^- +F^- G^+  \!\!&=&\!\! -\frac{1}{16z^{10}}
\! \left(\begin{array}{l}
  2(\BB_{-1}\frac{\pl Q_z}{\pl s})(\BB_{-1}^3Q_z )
    \\ \\
    -32(\BB_0-2s\frac{\pl}{\pl s}-2)\BB_{-1}^2Q_z
     \\ \\ +
      (\BB_{-1}^2\frac{\pl Q_z}{\pl s})(\BB_{-1}^2Q_z )-16 \frac{\pl ^3Q_z}{\pl s^3}
     \end{array}
      \right)
     +O(\frac{1}{z^{8}})
\\
   F^+\tilde\BB_{-1}^2F^-
-F^-\tilde\BB_{-1}^2F^+  \!\!&=&\!
 \vr\frac{\frac{\pl}{\pl
s}\BB_{-1}^3 Q_z}{4z^{7}}+O(\frac{1}{z^{5}})
  .\eean}
Using this expressions, one easily deduces for small
$z$,
{\footnotesize\bean
 0&=&\!\!\!\!\!\!\!\left.\left\{\!\!\begin{array}{l}
\Bigl(F^+\tilde\BB_{-1}G^-+F^-\tilde\BB_{-1}G^+ \Bigr)
\Bigl(F^+\tilde\BB_{-1}F^- -F^-\tilde\BB_{-1}F^+ \Bigr) \\
 -
 \Bigl(F^+ G^- +F^- G^+
\Bigr) \Bigl(F^+\tilde\BB_{-1}^2F^-
-F^-\tilde\BB_{-1}^2F^+ \Bigr)
\end{array}\!\!\!\right\}
 \right|_{\tiny\begin{array}{l}
                    u\mapsto \frac{\sqrt{2}}{z^2}
                    \sqrt{\frac{\frac{1}{2}+sz^2}{\frac{1}{2}-sz^2}} \\
                    v_i\mapsto x_i\frac{z\sqrt{2}}{\sqrt{\frac{1}{4}-s^2z^4}}
           \\
           n\mapsto \frac{2}{z^4}\end{array} }
\\  &&\\
 &=& -\frac{ \vr}{2 z^{17}}
 \left(\begin{array}{l}
   \left\{\BB_{-1}^2\frac{\pl Q}{\pl
s},\frac{1}{2}
   \frac{\pl^3 Q_z}{\pl s^3}+
(\BB_0-2)\BB_{-1}^2Q_z\right\}_{\BB_{-1}}
  \\  \\
  +\frac{1}{16}\BB_{-1} \frac{\pl Q_z}{\pl
s}\left\{\BB_{-1}^3Q_z,\BB_{-1}^2\frac{\pl Q_z}{\pl s}
\right\}_{\BB_{-1}}
  \end{array}\right)+O(\frac{1}{z^{15}})\\  &&\\
 &=& -\frac{ \vr}{2 z^{17}}
 \left(
   \begin{array}{l}
   \mbox{the same expression for
    $Q(s;x_1,\ldots, x_{2r})$}
\end{array}
  \right)+O(\frac{1}{z^{16}})
  ,\eean}
using (\ref{5.11}) in the last equality. Taking the
limit when $z\rightarrow 0$ yields equation (\ref{0.5})
of Theorem 0.2.\qed


\section{Appendix 1}

Setting
$$\mu_{i+j-1}(\pm
a):=\mu_{ij}^{\pm}(t,s,u;\beta,\BR)\big\vert_{t=s=u=\beta
=0}=\int_{\BR} z^{i+j-1}e^{-\frac{z^2}{2}\pm az}dz,
$$
one computes\footnote{Remember $n=k_1+k_2$.}:
\begin{lemma}
\bean
\tau_{k_1k_2}(t,s,u;\beta,\BR)\big\vert_{t=s=u=\beta
=0}&=& \det\left(\begin{array}{c}
\left(\displaystyle{\mu_{i+j}(a)}\right)_{\tiny{\begin{array}{l}
      0\leq i\leq k_1-1\\
      0\leq j\leq n-1
      \end{array}}}\\
\\
\left(\displaystyle{\mu_{i+j}(-a)}\right)_{\tiny{\begin{array}{l}
      0\leq i\leq k_2-1\\
      0\leq j\leq n-1
      \end{array}}}
\end{array}\right)\\
\\
&=&c_{k_1k_2}a^{k_1k_2}e^{\frac{(k_1+k_2)}{2}} a^2.
\eean with
$$
c_{k_1k_2}=(-2)^{k_1k_2}(2\pi)^{\frac{k_1+k_2}{2}}
\prod_0^{k_1-1}j!   \prod_0^{k_2-1}j!~.$$
\end{lemma}

\proof By explicit integration, one computes

$$
\mu_0(a)=\sqrt{2\pi}~e^{\frac{a^2}{2}}\mbox{~and~}\mu_i(\pm
a)=\sqrt{2\pi}\left(\pm\frac{d}{da}\right)^ie^{\frac{a^2}{2}}.
$$
Define the Hermite polynomials (except for a minor
change of variables)
$$
p_i(a):=e^{-\frac{a^2}{2}}\left(\frac{d}{da}\right)^ie^{\frac{a^2}{2}}=
\left(\frac{d}{da}+a\right)p_{i-1}(a).
$$
The following holds:
$$p_{2i}(a)=\mbox{even polynomial,}\qquad p_{2i+1}(a)=\mbox{odd
polynomial of $a$},
$$
which is used in equality $\stackrel{\ast\ast}{=}$
below, and
$$p_{k+n}(a)=p_k^{(n)}+\beta_1(a)p_k^{(n-1)}+\beta_2(a)p_k^{(n-2)}+\ldots
+\beta_np_k,
$$
where $p_k^{(n)}:=\left(\frac{d}{da}\right)^np_k$ and
where $\beta_i(a)$ are polynomials in $a$, independent
of $k$; this feature is used in equality
$\stackrel{\ast}{=}$ below. Then we compute:
%
%
\bean
\lefteqn{\tau_{k_1k_2}(t,s,u;\beta,\BR)\big\vert_{t=s=u=\beta =0}}\\
&=&(\sqrt{2\pi})^n~e^{\frac{na^2}{2}}\det\left(\begin{array}{c}
\left(p_{i+j}\right)_{\tiny{\begin{array}{l}
      0\leq i\leq k_1-1\\
      0\leq j\leq n-1
      \end{array}}}\\
\\
\left((-1)^{i+j}p_{i+j}\right)_{\tiny{\begin{array}{l}
      0\leq i\leq k_2-1\\
      0\leq j\leq n-1
      \end{array}}}
\end{array}\right)\\
\\
&=&(\sqrt{2\pi})^n(-1)^{\frac{k_2(k_2-1)}{2}}
e^{\frac{na^2}{2}}\det\left(\begin{array}{c}
\left(p_{i+j}\right)_{\tiny{\begin{array}{l}
      0\leq i\leq k_1-1\\
      0\leq j\leq n-1
      \end{array}}}\\
\\
\left((-1)^{j}p_{i+j}\right)_{\tiny{\begin{array}{l}
      0\leq i\leq k_2-1\\
      0\leq j\leq n-1
      \end{array}}}
\end{array}\right)\\
\\
&\stackrel{\ast}{=}&(\sqrt{2\pi})^n
(-1)^{\frac{k_2(k_2-1)}{2}}e^{\frac{na^2}{2}}\det
 \left(\begin{array}{c}
 \left(p_{j}^{(i)}\right)_{\tiny{\begin{array}{l}
      0\leq i\leq k_1-1\\
      0\leq j\leq n-1
      \end{array}}}\\
\\
\left((-1)^{j}p_{j}^{(i)}\right)_{\tiny{\begin{array}{l}
      0\leq i\leq k_2-1\\
      0\leq j\leq n-1
      \end{array}}}
\end{array}\right)\\
\\
&\stackrel{\ast\ast}{=}&c_{k_1k_2}e^{\frac{na^2}{2}} \
\det\left(
 \begin{array}{c}
 \left((a^{j-1})^{(i)}\right)_{\tiny{\begin{array}{l}
      0\leq i\leq k_1-1\\
      1\leq j\leq n
      \end{array}}}\\
\\
\left(((-a)^{j-1})^{(i)}\right)_{\tiny{\begin{array}{l}
      0\leq i\leq k_2-1\\
      1\leq j\leq n
      \end{array}}}
\end{array}\right)\\
\\
&=&c_{k_1k_2}e^{\frac{na^2}{2}}\det\left(\begin{array}{c}
\left(\al_{ij}a^{j-i}\right)_{\tiny{\begin{array}{l}
      1\leq i\leq k_1\\
      1\leq j\leq n
      \end{array}}}\\
\\
\left(\al_{ij}a^{j-i+k_1}\right)_{\tiny{\begin{array}{l}
      k_1+1\leq i\leq n\\
      1\leq j\leq n
      \end{array}}}
\end{array}\right)\\
\\
&=&c_{k_1k_2}e^{\frac{na^2}{2}}\sum_{\sg\in
S_n}(-1)^{\sg}\prod_{1\leq i\leq
k_1}\al_{i\sg(i)}a^{\sg(i)-i}\prod_{k_1+1\leq i\leq
n}\al_{i\sg(i)}a^{\sg(i)-i+k_1}\\
\\
&=&c_{k_1k_2}e^{\frac{na^2}{2}}\sum
(-1)^{\sg}a^{\sum_1^n(\sg(i)-i)}
 \left(a^{k_1}\right)^{k_2}
\prod_{1\leq i\leq n}\al_{i\sg(i)}\\
\\
&=&c'_{k_1k_2}e^{\frac{(k_1+k_2)}{2}a^2}a^{k_1k_2},
\eean where the $\al_{ij}$ are coefficients, some of
which vanish. Indeed, each of the blocks in the matrix
above is upper-triangular. To evaluate $c'_{k_1k_2}$,
observe, upon completing the squares in the exponentials
and setting $x_j\mapsto x_j-a,~y_j\mapsto y_j+a$ in the
integral,
 \bean
  \lefteqn{  \tau_{k_1k_2}(t,s,u;\beta;\BR)
  \Biggr|_{t=s=u=\beta=0}}\\
   &=&\frac{1}{k_1!k_2!}
  \int_{\BR^{k_1+k_2 }}\Dt_{k_1+k_2 }(x,y)
\left(\Dt_{k_1}(x)\prod^{k_1}_{j=1}
  e^{-\frac{x_j^2}{2}+ax_j }
dx_j\right)
 \no\\ &&~~~
~~~~~~~~~~~~~~~~~~~~~~~~~~~~~~~~~\left(\Dt_{k_2}(y)
\prod^{k_2}_{j=1}
e^{-\frac{y_j^2}{2}-ay_j}
 dy_j\right)
  .~~~~~
  \eean
This integral equals
 \bean
  &=&
   e^{(k_1+k_2)a^2/2}\left(
(-2a)^{k_1k_2}c_{k_1,0}c_{0,k_2}+\mbox{lower order terms
in $a$} \right) \\
 &=&e^{(k_1+k_2)a^2/2}
\left((-2a)^{k_1k_2}(2\pi)^{\frac{k_1+k_2}{2}}
\prod_0^{k_1-1}j!\prod_0^{k_2-1}j!
 +\mbox{lower order terms
in $a$} \right)
   \eean
The result in the first part of this proof implies the
absence of the lower terms and thus Lemma 6.1.\qed


\section{Appendix 2}

In order to compute the asymptotics (\ref{5.10}) for the
expression $G^{\vr}$ and $\BB_{-1}G^{\vr}$, as defined
in (\ref{0.7}), one needs the following asymptotics:
{\footnotesize\bean
  F^{\vr}& =& -\frac{4}{z^4}-\frac{1}{4z^2}\BB_{-1}^2 Q_z+
 \frac{\vr}{4  z}\BB_{-1}\frac{\pl  Q_z}{ \pl
 s }+O(z)
 \\  \\
  \frac{1}{\sqrt{2}}\tilde\BB_{-1}F^{\vr}
 &=&
 -\frac{1}{16z^3}\BB_{-1}^3 Q_z
 +
\frac{\vr}{16  z^2}\BB_{-1}^2\frac{\pl Q_z}{\pl s}
 -\frac{\vr s}{8   }\BB_{-1}^2\frac{\pl Q_z}{\pl s}
   +  O(z)
 \\  \\
   \tilde\BB_{-1}^2F^{\vr}
 &=&
 -\frac{1}{32z^4}\BB_{-1}^4 Q_z
 +
\frac{\vr}{32  z^3}\BB_{-1}^3\frac{\pl Q_z}{\pl s}
 -\frac{\vr s}{16  z }\BB_{-1}^3\frac{\pl Q_z}{\pl s}
   +  O(1)
  \no\\
  \frac{1}  {\sqrt{2}}\frac{\pl}{\pl a}F^{\vr}&=&
   -\frac{\vr}{16z^2}\BB_{-1}^2 \frac{\pl  Q_z}{ \pl
 s }
  + O(\frac{1}{z})
  \\ \\
   \frac{\pl}{\pl a}\tilde\BB_{-1}
  F^{\vr}&=&
   -\frac{\vr}{32z^3}\BB_{-1}^3 \frac{\pl  Q_z}{ \pl
 s }
 +O(\frac{1}{z^2 })
  \\ \\
 \frac{1}{\sqrt{2}}H_1^{\vr}& =&  \frac{6\vr}{z^6}
 +\frac{4\vr s}{z^4}
  - \frac{1}{8  z^3}\BB_{-1}\frac{\pl  Q_z}{ \pl
 s }+O(\frac{1}{z^2})
 \eean\bean
  \tilde\BB_{-1}H_1^{\vr}
 &=&
 -\frac{1}{16z^4}\BB_{-1}^2 \frac{\pl  Q_z}{ \pl
 s }
 -
\frac{\vr}{16  z^3}\BB_{-1} \frac{\pl^2 Q_z}{\pl s^2}
 +\frac{1}{8z^2   }\BB_0\BB_{-1}^2Q_z
   +  O(\frac{1}{z })
 \\  \\
 \frac{1}{\sqrt{2}} \tilde\BB_{-1}^2H_1^{\vr}
 &=&
 -\frac{1}{64z^5}\BB_{-1}^3 \frac{\pl  Q_z}{ \pl
 s }
 -
\frac{\vr}{64  z^4}\BB_{-1}^2\frac{\pl^2 Q_z}{\pl s^2}
 +\frac{ 1}{32  z^3}(\BB_{0}+1)\BB_{-1}^3Q_z
   +  O(\frac{1}{z^2})
  \\  \\
   \frac{1}{\sqrt{2}}H_2^{\vr}& =&
   \frac{ \vr}{4  z^4}\BB_{-1}^2 Q_z
  +    O(\frac{1}{z^3})
   \\  \\
    \frac{\pl}{\pl a}H_2^{\vr}&=&
   \frac{1}{8z^4}\BB_{-1}^2 \frac{\pl  Q_z}{ \pl
 s }
  -\frac{1}{16z^3}\BB_{-1}  \frac{\pl^2  Q_z}{ \pl
 s^2 }-\frac{1}{16z^2}(\frac{\pl^3}{\pl s^3}
  -4\BB_{-1}^2 )Q_z+ O(\frac{1}{z  })
 \\  \\
  \tilde\BB_{-1}H_2^{\vr}
 &=&
 \frac{\vr}{8z^5}\BB_{-1}^3 Q_z
   + O(\frac{1}{z^4 })
  \\ \\
   \frac{1}{\sqrt{2}} \frac{\pl}{\pl a}\tilde\BB_{-1}
  H_2^{\vr}&=&
    \frac{1}{32z^5}\BB_{-1}^3 \frac{\pl  Q_z}{ \pl
 s }
   -\frac{\vr}{64z^4}\BB_{-1}^2\frac{\pl^2  Q_z}{ \pl
 s^2 }
  -\frac{1}{64z^3}(  \frac{\pl^3   }{ \pl
 s^3 }-4\BB_{-1}^2)\BB_{-1}Q_z+O(\frac{1}{z ^2})
   .\eean}
%

\end{document}